\documentclass{article}
\usepackage{amsfonts}
\usepackage{amsmath}

\newtheorem{theorem}{Theorem}

\newtheorem{corollary}[theorem]{Corollary}

\newtheorem{definition}[theorem]{Definition}
\newtheorem{example}[theorem]{Example}

\newtheorem{lemma}[theorem]{Lemma}

\begin{document}

\title{ Polling systems and multitype branching processes in random
environment counted by random characteristics }
\author{V.A. Vatutin\thanks{%
Research supported by the grant RFBR 08-01-00078} \\
Steklov Mathematical Institute, Gubkin street, 8, 119991,\\
Moscow,Russia\\
E-mail: \texttt{vatutin@mi.ras.ru}}
\maketitle

\begin{abstract}
By the methods of multitype branching processes in random environment
counted by random characteristics we study the tail distribution of busy
periods and some other characteristics of the branching type polling systems
in which the service disciplines, input parameters and service time
distributions are changing in a random manner.

\end{abstract}

\section{\protect\bigskip Polling systems with service policies of branching
type\label{PollBranch}}

We consider a polling system consisting of a single server and $m$
stations with infinite-buffer queues indexed by $i\in \left\{
1,\ldots ,m\right\}.$ Initially there are no customers in the
system. When customers arrive to the system the server starts
immediately the service by visiting the stations in cyclic order
$(1\rightarrow 2\rightarrow \cdots \rightarrow m\rightarrow
1\rightarrow \cdots )$ starting at station $1$ according to a
selected service policy (to be described later on) and with zero
switchover times between queues. Later on the initial stage of
services $(1\rightarrow 2\rightarrow \cdots \rightarrow m)$ will
be called the zero cycle. The subsequent routes of the server will
be called the first cycle, the second cycle and so on. When the
system is empty the server waits their arrival at a parking place
$R$. Customers arrive to the queues in accordance with a point
process whose parameters are changing in a random manner each time
when the server switches from station to station .

To give a rigorous description of the arrival and service processes for the
system in question we need some  notions. Let $\mathbf{s:}%
=(s_{1},\ldots ,s_{m})\in \left[ 0,1\right] ^{m}$ be a $m$-dimensional
variable,
\begin{equation*}
\mathbf{s}^{\mathbf{k}}:=s_{1}^{k_{1}}\cdots s_{m}^{k_{m}},\,k_{i}\in
\mathbb{N}_{0}=\left\{ 0,1,2,\ldots \right\} ,
\end{equation*}%
$\mathbf{e}_{i}:=\left( 0,\ldots ,1,\ldots ,0\right) $ be a $m-$dimensional
vector with zero components except the $i$-th equal to $1,$ $\mathbf{0}$ and
$\mathbf{1}$ are $m-$dimensional vectors whose all components are equal to $%
0 $ and $1$, respectively. Let
\begin{equation*}
\phi ^{(i)}(\mathbf{s;}\lambda ):=\mathbf{E}\left[ s_{1}^{\theta
_{i1}}s_{2}^{\theta _{i2}}\cdots s_{m}^{\theta _{im}}e^{-\lambda \phi _{i}}%
\right] ,i=1,...,m
\end{equation*}%
be the mixed probability generating function (m.p.g.f.'s) of $(m+1)-$
dimensional vectors $\left( \theta _{i1},\ldots ,\theta _{im};\phi
_{i}\right) $, where $\theta _{ij}$ are nonnegative integer-valued random
variables and $\phi _{i}$ is a nonnegative random variable. Denote%
\begin{equation*}
\mathbf{\phi }(\mathbf{s;}\lambda ):=\left( \phi ^{(1)}(\mathbf{s;}%
\lambda ),\ldots ,\phi ^{(m)}(\mathbf{s;}\lambda )\right)
\end{equation*}%
the respective vector-valued m.p.g.f.\ and let
$\mathcal{F}=\left\{ \mathbf{\phi }(\mathbf{s;}\lambda )\right\} $
be the set of all such m.p.g.f. Assume that a probability
measure $\mathbb{P}$ is specified on the natural $\sigma $-algebra of $%
\mathcal{F}$  and let
\begin{equation*}
\mathbf{\phi }_{0}(\mathbf{s;}\lambda ),\mathbf{\phi }_{1}(\mathbf{s;}%
\lambda ),\mathbf{\phi }_{2}(\mathbf{s;}\lambda ),\ldots ,
\end{equation*}%
be a sequence selected in iid manner from $\mathcal{F}$ in accordance with $%
\mathbb{P},$ where $\mathbf{\phi }_{n}(\mathbf{s;}\lambda ):=\left( \phi
_{n}^{(1)}(\mathbf{s;}\lambda ),\ldots ,\phi _{n}^{(m)}(\mathbf{s;}\lambda
)\right) $ with%
\begin{equation*}
\phi _{n}^{(i)}(\mathbf{s;}\lambda ):=\mathbf{E}\left[ s_{1}^{\theta
_{i1}(n)}s_{2}^{\theta _{i2}(n)}\cdots s_{m}^{\theta _{im}(n)}e^{-\lambda
\phi _{i}(n)}\right] .
\end{equation*}

In the present paper we investigate such polling systems whose
arrival and service procedures of customers meet the following
property.

\textbf{Branching property}. When the server arrives to station $i$ for the $%
n$-th time and find there, say, $k_{i}$ customers labelled $1,2,...,k_{i}$ ,
then, during the course of the server's visit, the arrival of customers is
arranged in such a way that after the end of each stage of service of
customer $j$ (the number of such stages may be more than one if the service
discipline admits feedback) the queues at the system will be increased by a
random population of customers $\left( \theta _{i1}(n,j),\ldots ,\theta
_{im}(n,j)\right) ,$ where $\theta _{il}(n,j)$ is the number of customers
added to the $l$-th station, and, in addition, a final product of size $\phi
_{i}(n,j)\geq 0$ will be added to the system. It is assumed that the vectors
$\left( \theta _{i1}(n,j),\ldots ,\theta _{im}(n,j);\phi _{i}(n,j)\right) ,$
$j=1,2,\ldots ,k_{i}$ are iid and such that%
\begin{equation*}
\mathbf{E}\left[ s_{1}^{\theta _{i1}(n,j)}s_{2}^{\theta _{i2}(n,j)}\cdots
s_{m}^{\theta _{im}(n,j)}e^{-\lambda \phi _{i}(n,j)}\right] =\phi _{n}^{(i)}(%
\mathbf{s;}\lambda ).
\end{equation*}

Note that $\phi _{n}^{(i)}(\mathbf{s;}\lambda )$ is a random m.p.g.f. and,
therefore, the parameters of the polling system are changed from cycle to
cycle in a random manner. Since in the paper we are interesting in the
characteristics related with busy periods of the system whose work starts by
the arrival of a single customer to a station $J\in \left\{ 1,...,m\right\} $
at moment 0, the law of arrival of customers to an idle system plays no role
for subsequent arguments.

The aim of the article is to study the distribution of the total size of the
final product accumulated in the system during a busy period of the server.
In particular, letting the final product $\phi _{i}(n,j)$ be the service
time of the $j$-th customer served during the $n$-th visit of the server to
station $i$, we provide conditions under which the tail distribution of the
length of the busy period decays, as $y\rightarrow \infty ,$ like $%
const\times y^{-\kappa }$ for some $\kappa >0$.

Before we proceed to the rigorous statements of our results,
consider two examples of branching type polling systems.

Let $\mathcal{T}_{+}\mathcal{=}\left\{ T\right\} $ be the set of
all probability distributions of nonnegative random variables,  $
\mathcal{T}_{+}^{m}:=\left\{ \left( T_{1},...,T_{m}\right)
:T_{i}\in \mathcal{T}_{+}\right\} $ be the set of all
$m$-dimensional tuples of such distributions,
$\mathcal{M}_{\varepsilon }=\left\{ \mathcal{E}\right\} $ be the
set of all $m\times m$ matrices $\mathcal{E}=\left( \varepsilon
_{ij}\right) _{i,j=1}^{m}$ with nonnegative elements, and
$\mathcal{M}_{\gamma }=\left\{ \Gamma \right\} $ be the set of all
$m\times (m+1)$ matrices $\Gamma =\left( \gamma _{ij}\right)
_{i=1,j=0}^{m}$ with nonnegative elements such that

\begin{equation*}
\sum_{j=0}^{m}\gamma _{ij}=1,\quad i=1,...,m.
\end{equation*}

Let $P$\ be a measure on the Borel $\sigma -$algebra of the space $\mathcal{M%
}_{\varepsilon }\times \mathcal{M}_{\gamma }\times \mathcal{T}_{+}^{m}.$

\begin{example}
\label{Examp1}(motivated by \cite{PMPP2008}). Consider a polling
system with $m$ stations and a single server performing cyclic
service of the customers at the stations. Assume that initially
there are no customers in the system and the server is located at
parking place $R$. Assume that given the idle system the flow of
customers arriving to station $i$ is Poisson with, say, a
deterministic rate $\varepsilon _{i}$. When the first customer
appears in the system the server selects a random element $\left(
\mathcal{E}_{0},\Gamma
_{0},\mathbf{T}_{0}\right) \in \mathcal{M}_{\varepsilon }\times \mathcal{M}%
_{\gamma }\times \mathcal{T}_{+}^{m}$ with%
\begin{equation*}
\mathcal{E}_{0}\mathcal{=}\left( \varepsilon _{ij}(0)\right)
_{i,j=1}^{m},\quad \Gamma_{0}=\left( \gamma _{ij}(0)\right)
_{i=1,j=0}^{m},\quad \mathbf{T}_{0}=\left( T_{10},...,T_{m0}\right)
\end{equation*}%
and immediately starts its zero service cycle $(1\rightarrow
2\rightarrow \cdots \rightarrow m)$ adopting the gated server
policy with zero switchover times. Namely, the server serves all
the customers that were queueing at a station when the server
arrived and then instantly jumps to the next (in cyclic order)
station. For the period while the server performs the batch of
services at station $i$, new customers arrive to the system
according to independent Poisson flows with intensities given by
the vector $(\varepsilon _{i1}(0),\varepsilon
_{i2}(0),...,\varepsilon _{im}(0))$ (some the components may be
equal to zero) and the service
times of customers are iid and distributed according to $T_{i0}(x):=\mathbf{P%
}\left( \tau _{i}(0)\leq x\right) $. Each served customer either goes to
station $j\in \{1,\ldots ,m\}$ with probability $\gamma _{ij}(0)$ or leaves
the system with probability $\gamma _{i0}(0)$ independently of other events.
Besides, after the end of each service period of a customer the customer
contributes\ to the system its service time as the final product. It is
assumed that given $\left( \mathcal{E}_{0},\Gamma _{0},\mathbf{T}_{0}\right)
$ the service times and the arrival process of new customers are independent.

The subsequent routes $n=1,2,\ldots $ have the same probabilistic structure
specified by the tuples
\begin{equation*}
\mathcal{E}_{n}\mathcal{=}\left( \varepsilon _{ij}(n)\right)
_{i,j=1}^{m},\quad \Gamma_{n}=\left( \gamma _{ij}(n)\right)
_{i=1,j=0}^{m},\quad \mathbf{T}_{n}=\left( T_{1n},...,T_{mn}\right)
\end{equation*}%
with only difference that at the beginning of cycle $n\geq 1$ there is a
possibility to have more than one customer in the system.
\end{example}

Let us show that this system possesses a branching property. To this aim
denote%
\begin{equation*}
t_{in}(\lambda ):=\int_{0}^{\infty }e^{-\lambda x}dT_{in}(x)
\end{equation*}%
the Laplace transform of the distribution $T_{in}(x)$ of the random variable
$\tau _{i}(n)$ representing the service time of a customer at station $i$
during the $n-$th visit of the server.

It is not difficult to check that the m.p.g.f. $\mathbf{\phi }_{n}(\mathbf{s;%
}\lambda )$ has the components (in our setting and the service time of
customers as the final product)
\begin{eqnarray}
\phi _{n}^{(i)}(\mathbf{s;}\lambda ) &:&=\mathbf{E}\left[ s_{1}^{\theta
_{i1}(n)}s_{2}^{\theta _{i2}(n)}\cdots s_{m}^{\theta _{im}(n)}e^{-\lambda
\tau _{i}(n)}\right]  \notag \\
&=&\int_{0}^{\infty }\mathbf{E}\left[ s_{1}^{\theta _{i1}(n)}s_{2}^{\theta
_{i2}(n)}\cdots s_{m}^{\theta _{im}(n)}|\tau _{i}(n)=x\right] e^{-\lambda
x}dT_{in}(x)  \notag \\
&=&\left( \gamma _{i0}(n)+\sum_{j=1}^{m}\gamma _{ij}(n)s_{j}\right)
\int_{0}^{\infty }\prod_{j=1}^{m}e^{\varepsilon
_{ij}(n)(s_{j}-1)x}e^{-\lambda x}dT_{in}(x)  \notag \\
&=&\left( \gamma _{i0}(n)+\sum_{j=1}^{m}\gamma _{ij}(n)s_{j}\right)
t_{in}\left( \lambda +\sum_{j=1}^{m}\varepsilon _{ij}(n)(1-s_{j})\right)
\label{GateTotal}
\end{eqnarray}%
and%
\begin{equation}
h_{n}^{(i)}(\mathbf{s}):=\phi _{n}^{(i)}(\mathbf{s;}0)=\left( \gamma
_{i0}(n)+\sum_{j=1}^{m}\gamma _{ij}(n)s_{j}\right) t_{in}\left(
\sum_{j=1}^{m}\varepsilon _{ij}(n)(1-s_{j})\right) .  \label{GateSimple}
\end{equation}

\begin{example}
\label{Examp2}(compare with \cite{PMPP2008}). Consider the same polling
system as earlier but assume now that at each station the server adopts the
exhaustive server policy: it serves all the customers that were queueing at
the station when the server arrived together with all subsequent arrivals up
until the queue becomes empty and then instantly jumps to the next station.
\end{example}

Let us show that this system possesses a branching property as well.

Since during the $n$-th visit of the server to station $i$ each customer at
this station is served by the server a random number of times having shifted
by~1 geometric distribution with parameter $\gamma _{ii}(n)$, the Laplace
transform $w_{in}(\lambda )$ of the distribution of the random variable $%
\eta _{i}(n)$, the total service time of a customer at station $i$ during
the $n-$th cycle, has the form
\begin{equation}
w_{in}(\lambda )=\frac{1-\gamma _{ii}(n)}{1-\gamma _{ii}(n)t_{in}(\lambda )}%
t_{in}(\lambda ).  \label{ForW}
\end{equation}%
Let now $\left( \sigma _{i1}(n),\ldots ,\sigma _{im}(n)\right) $ be a random
vector distributed as the vector of the number of customers arriving to the
stations of the polling system during the total service time of a customer
at queue $i$ which, after the end of its service at station $i$ moves either
to a station $j\neq i$ or leaves the system. Denote%
\begin{equation*}
y_{in}(\mathbf{s}):=\frac{\gamma _{i0}(n)+\sum_{j\neq i}\gamma _{ij}(n)s_{j}%
}{1-\gamma _{ii}(n)}.
\end{equation*}%
Similarly to (\ref{GateTotal}) one can show that
\begin{equation}
\mathbf{E}\left[ s_{i}^{\sigma _{i1}(n)}\cdots s_{m}^{\sigma
_{im}(n)}e^{-\lambda \tau _{i}(n)}\right] =y_{in}(\mathbf{s})w_{in}\left(
\lambda +\sum_{j=1}^{m}\varepsilon _{ij}(n)(1-s_{j})\right) .  \label{Stand}
\end{equation}

Let now $l_{in}(\lambda )$ be the Laplace transform of the distribution of a
busy period $\eta _{i,tot}(n)$ of the system generated by a single customer
in an $M/G/1$ queue with arrival rate $\varepsilon _{i}(n)$ and the service
time distributed as $\eta _{i}(n),$ and let $\left( \theta _{i1}(n),\ldots
,\theta _{im}(n)\right) $ be the total number of new customers arriving to
the $i$-th stations of our polling system within the time-interval
distributed as $\eta _{i,tot}(n)$. In this case $l_{in}(\lambda )$ is a
unique solution of the equation%
\begin{equation*}
l_{in}(\lambda )=w_{in}\left( \lambda +\varepsilon _{ii}(n)\left(
1-l_{in}(\lambda )\right) \right)
\end{equation*}%
and (compare with \cite{Res93}) the functions%
\begin{equation*}
\phi _{n}^{(i)}(\mathbf{s;}\lambda ):=\mathbf{E}\left[ s_{1}^{\theta
_{i1}(n)}s_{2}^{\theta _{i2}(n)}\cdots s_{m}^{\theta _{im}(n)}e^{-\lambda
\eta _{i,tot}(n)}\right] ,i=1,2,\ldots ,m
\end{equation*}%
are (in our setting) unique solutions of the equations%
\begin{equation*}
\phi _{n}^{(i)}(\mathbf{s;}\lambda )=y_{in}(\mathbf{s})w_{in}\left( \lambda
+\sum_{j\neq i}\varepsilon _{ij}(n)(1-s_{j})+\varepsilon _{ii}(n)(1-\phi
_{n}^{(i)}(\mathbf{s;}\lambda ))\right) .
\end{equation*}%
In particular,%
\begin{equation}
h_{n}^{(i)}(\mathbf{s})=y_{in}(\mathbf{s})w_{in}\left( \sum_{j\neq
i}\varepsilon _{ij}(n)(1-s_{j})+\varepsilon _{ii}(n)(1-h_{n}^{(i)}(\mathbf{s}%
))\right) .  \label{Exh1}
\end{equation}

Note that in the framework of the suggested approach one can consider also
models with batch arrivals of customers as well. In this case one should
replace, for instance, everywhere in $t_{in}(\cdot )$ variables $%
s_{1},\ldots ,s_{m}$ by the respective probability generating functions of
the sizes of batches of customers arriving to stations $j=1,\ldots ,m$.

Since we imposed the branching property on the service disciplines
of the systems it is not a surprise for the reader that we conduct
our investigation by the methods of the theory of branching
processes. This approach is not new. Polling systems possessing
the branching property in which the probability generating
functions $h_{n}^{(i)}(\mathbf{s}),\quad n=1,2,\ldots$ are
nonrandom and the same for all cycles were considered in
particular, in \cite{Fur91}, \cite{Res93}, \cite{VD2002} and quite
recently in \cite{YY07}, \cite{Alt2009} and \cite{AltF2007}. These
models cover many classical service policies, including the
exhaustive, gated, binomial-gated and their feedback modifications
(see surveys   \cite{ViSe06} and \cite{Mei07} for definitions and
more details).

Polling systems with input parameters and service disciplines
changing in a random manner and(or) depending on the states of
systems are not studies yet in full generality. The analysis of
such systems uses rather often the fluid method
(\cite{FC98}-\cite{FL98}) or a method based on the construction of
appropriate Lyapunov functions (\cite{PMPV2006} - \cite{PMPP2008})
 (in the last case without
reference to branching processes). The authors of \cite{PMPP2008}
write that it is possible to generalize their models. However
"...This leads to a lot of complexities in the proofs." In the
present paper we show that the reduction of the problems related
to the branching type polling systems with final product and
evolving in random environment (BTPSFPRE) to the respective
problems for the multitype branching processes counted by random
characteristics and evolving in random environment (MBPRCRE) gives
the desired answers by a unique method and in a general situation.
Since m.p.g.f.'s are selected for each cycle independently of the
past, different service disciplines are allowed at different
stations in our model. Moreover, the service disciplines at the
stations may be changed at random from visit to visit. In
particular, we allow for the mixture of exhaustive and gated
service disciplines. Here we consider the models with zero
switchover times. The study of BTPSFPRE with positive (and random)
switchover times may be conducted by a similar method. This,
however, requires more efforts and will be done elsewhere.

The scheme of the remaining part of the paper looks as follows. Section \ref%
{SecBPRE} is devoted to a detailed description of multitype branching
processes counted by random characteristics and evolving in random
environment. In Section \ref{SecLimT} we recall some known results for
ordinary multitype branching processes in random environment (MBPRE) and
formulate the main results of our paper describing the asymptotic behavior
of the total size of the final product for subcritical MBPRCRE. Section \ref%
{SecAux} recalls an important statement related with the asymptotic
properties of the tail distribution of  infinite sums of products of random matrices. In Section %
\ref{NonRand} we deduce some estimates for the moments of the population
sizes of subcritical MBPRE's. The proofs of the main results of our paper
are collected in Section \ref{MMain}. And, finally, we demonstrate in
Section \ref{SecPoll} how one can use the obtained results for MBPRCRE to
make conclusions about the probabilistic properties of various
characteristics of the branching type polling systems.

\section{Branching processes in random environment counted by random
characteristics\label{SecBPRE}}

As we claimed in the previous section, the main goal of the present paper is
to analyze properties of busy periods and some other characteristics of the
branching type polling systems by means of MBPRCRE. However, to start such
analysis we need to pass through a relatively long way of notation and
statements which is practically always the case when one consider multitype
branching processes.

Let $\left( \xi _{1},\ldots ,\xi _{m};\varphi \right) $ be a $(m+1)$%
-dimensional vector where the components $\xi _{1},\ldots ,\xi _{m}$ are
integer-valued nonnegative random variables and $\varphi $ is a nonnegative
random variable and let%
\begin{equation*}
F(\mathbf{s;}\lambda )=\mathbf{E}\left[ s_{1}^{\xi _{1}}s_{2}^{\xi
_{2}}\cdots s_{m}^{\xi _{m}}e^{-\lambda \varphi }\right] ,\,\mathbf{s}%
=(s_{1},\ldots ,s_{m})\in \left[ 0,1\right] ^{m},\,\lambda \geq 0,
\end{equation*}%
be the respective m.p.g.f. Denote $\ \mathcal{F}_{\lambda }:=\left\{ F(%
\mathbf{s;}\lambda ),\mathbf{s}\in \left[ 0,1\right] ^{m}\right\} $ the set
of all such m.p.g.f.'s and let%
\begin{equation*}
\mathcal{F}_{\lambda }^{m}:=\mathcal{F}_{\lambda }\mathcal{\times \mathcal{F}%
_{\lambda }\times \cdots \times \mathcal{F}_{\lambda }=}\left\{ \mathbf{F}(%
\mathbf{s;}\lambda )=\left( F^{(1)}(\mathbf{s;}\lambda ),\ldots ,F^{(m)}(%
\mathbf{s;}\lambda )\right) \right\}
\end{equation*}%
be the $m$-times direct product of $\mathcal{F}_{\lambda }$. Let,
further,
\begin{equation*}
\mathcal{F}_{0}:=\left\{ f(\mathbf{s})=F(\mathbf{s;}0),\mathbf{s}%
\in \left[ 0,1\right] ^{m}\right\}
\end{equation*}%
be the set of all ordinary probability generating functions (p.g.f.'s)
\begin{equation*}
f(\mathbf{s})=\mathbf{E}\left[ s_{1}^{\xi _{1}}s_{2}^{\xi _{2}}\cdots
s_{m}^{\xi _{m}}\right]
\end{equation*}%
and
\begin{equation*}
\mathcal{F}^{m}_0:=\mathcal{F}_0\times \cdots \times \mathcal{F}_0=\left\{ \mathbf{f}%
(\mathbf{s})=\left( f^{(1)}(\mathbf{s}),\ldots
,f^{(m)}(\mathbf{s})\right) \right\}
\end{equation*}%
be the set of all $m$-dimensional (vector-valued) p.g.f.'s. Assume that a
probability measure $\mathbb{P}$ is specified on the natural $\sigma $%
-algebra $\mathcal{A}$ generated by the subsets of $\mathcal{F}_{\lambda
}^{m}.$ Let%
\begin{equation*}
\mathbf{F}_{0}(\mathbf{s;}\lambda ),\mathbf{F}_{1}(\mathbf{s;}\lambda
),\ldots ,\mathbf{F}_{k}(\mathbf{s;}\lambda ),\ldots \text{ with }\mathbf{F}%
_{n}(\mathbf{s;}\lambda ):=\left( F_{n}^{(1)}(\mathbf{s;}\lambda ),\ldots
,F_{n}^{(m)}(\mathbf{s;}\lambda )\right)
\end{equation*}%
be a sequence of vector-valued m.p.g.f.'s selected from $\mathcal{F}%
_{\lambda }^{m}$ in an iid manner in accordance with measure $\mathbb{P}$ .
The sequence $\left\{ \mathbf{F}_{n}(\mathbf{s;}\lambda ),n\geq 0\right\} $
is called a \textit{random environment}. \ With m.p.g.f. $F_{n}^{(i)}(%
\mathbf{s;}\lambda )$ we associate a random vector of offsprings $\mathbf{%
\xi }_{i}(n)=:\left( \xi _{i1}(n),\xi _{i2}(n),\ldots ,\xi _{im}(n)\right) $
and a random variable $\varphi _{i}(n)$ such that%
\begin{equation*}
F_{n}^{(i)}(\mathbf{s;}\lambda )=\mathbf{E}\left[ s_{1}^{\xi
_{i1}(n)}s_{2}^{\xi _{i2}(n)}\cdots s_{m}^{\xi _{im}(n)}e^{-\lambda \varphi
_{i}(n)}\right] .
\end{equation*}%
Now we may give an informal description of the MBPRCRE
\begin{equation*}
\mathbf{R}(n)=\left( \mathbf{Z}(n);\Phi (n)\right) ,n=0,1,\ldots ,
\end{equation*}%
which may be treated as the process describing the evolution of a $m-$type
population of particles with accumulation of a final product.

The starting conditions of the process are: \ a vector (may be random) of
particles $\mathbf{Z}(0)=\left( Z_{1}(0),\ldots ,Z_{m}(0)\right) $ where $%
Z_{i}(0)$ denotes the number of particles of type $i\in\{1,\ldots ,m\}$ in
the process at moment $0$, and an amount $\Phi (0)$ (may be random) of a
final product. All the particles have the unit life length and just before
the death produce children and final products independently of each other.
For instance, a particle, say, of type $i$ produces particles of different
types and adds some amount of the final product to the existing amount of
the final product in accordance with m.p.g.f. $F_{0}^{(i)}(\mathbf{s;}%
\lambda )$. The newborn particles constitute the first generation of the
MBPRCRE, have the unit life-length and dying produce, independently of each
other, offsprings and final products in accordance with their types and
subject to the m.p.g.f. $F_{1}^{(i)}(\mathbf{s;}\lambda ),i=1,2,\ldots ,m$
and so on.

A rigorous definition of the process we are interesting in looks as follows.

\begin{definition}
A $\left( m+1\right)$-dimensional Galton-Watson branching process
\begin{equation*}
\mathbf{R}_{\varphi }(n)=\mathbf{R}(n):=\left( \mathbf{Z}(n);\Phi (n)\right)
=\left( Z_{1}(n),\ldots ,Z_{m}(n);\Phi (n)\right) ,\,n\in \mathbb{N}_{0},
\end{equation*}%
counted by a random characteristics $\varphi $ in a fixed (but picked at
random) environment $\left\{ \mathbf{F}_{n}(\mathbf{s;}\lambda ),n\geq
0\right\} $ is a time-inhomogeneous Markov process with the state space
\begin{equation*}
\mathbb{N}_{0}^{m}\times \mathbb{R}_{0}:=\left\{ z=(z_{1},\ldots
,z_{m};w),\,z_{i}\in \mathbb{N}_{0};w\in \lbrack 0,\infty )\right\}
\end{equation*}%
defined as%
\begin{equation*}
\mathbf{R}(0)=\left( \mathbf{Z}(0);\Phi (0)\right) =(\mathbf{z};\varphi_0 )%
\mathbf{,}
\end{equation*}%
\begin{equation}
\;\mathbf{E}\left[ \mathbf{s}^{\mathbf{Z}(n+1)}e^{-\lambda \Phi (n+1)}|\;%
\mathbf{F}_{0},\ldots ,\mathbf{F}_{n};\mathbf{R}(0),\ldots ,\mathbf{R}(n)%
\right] =e^{-\lambda \Phi (n)}\left( \mathbf{F}_{n}(\mathbf{s;}\lambda
)\right) ^{\mathbf{Z}(n)}.
\end{equation}
\end{definition}

Note that the initial value $\mathbf{R}(0)$ may be random, and, for the
reason of applications to queueing systems we do not exclude the case $%
\mathbf{z}=\mathbf{0}$.

In what follows, to simplify notation, we write
\begin{equation*}
\mathbf{E}_{\mathbf{F}}\left[ \mathbf{s}^{\mathbf{Z}(n+1)}e^{-\lambda \Phi
(n+1)}\right] :=\mathbf{E}\left[ \mathbf{s}^{\mathbf{Z}(n+1)}e^{-\lambda
\Phi (n+1)}|\;\mathbf{F}_{0},\ldots ,\mathbf{F}_{n}\right] .
\end{equation*}%
Thus, $\mathbf{R}(0)=(\mathbf{z};\varphi_0),$
\begin{equation*}
\;\mathbf{E}_{\mathbf{F}}\left[ \mathbf{s}^{\mathbf{Z}(n+1)}e^{-\lambda \Phi
(n+1)}|\;\mathbf{R}(0),,\ldots ,\mathbf{R}(n)\right] =e^{-\lambda \Phi
(n)}\left( \mathbf{F}_{n}(\mathbf{s;}\lambda )\right) ^{\mathbf{Z}(n)},
\end{equation*}%
or%
\begin{eqnarray}
\mathbf{R}(n+1) &=&\left( \mathbf{0};\Phi (n)\right)
+\sum_{i=1}^{m}\sum_{k=1}^{Z_{i}(n)}\left( \mathbf{\xi }_{i}(n;k);\varphi
_{i}(n;k)\right) ,  \notag \\
\Phi (n+1) &=&\Phi
(0)+\sum_{l=0}^{n}\sum_{i=1}^{m}\sum_{k=1}^{Z_{i}(l)}\varphi _{i}(l;k)
\label{DDphi}
\end{eqnarray}%
where the random vector $\left( \mathbf{\xi }_{i}(n;k);\varphi
_{i}(n;k)\right)
$ represents the offspring vector and the size of the final product of the $%
k-$th particle of type $i$ of the $n-$th generation of the
process. Given the environment and $n=0,1,\ldots,$ and
$i\in\{1,\ldots,m\}$ the vectors
\begin{equation*}
\left( \mathbf{\xi }_{i}(n;k);\varphi _{i}(n;k)\right) ,k=1,2,\ldots
,Z_{i}(n)
\end{equation*}%
are independent and identically distributed: $\left( \mathbf{\xi }%
_{i}(n;k);\varphi _{i}(n;k)\right) \overset{d}{=}\left( \mathbf{\xi }%
_{i}(n);\varphi _{i}(n)\right) $.

Observe that if $\Phi (0)=0$ and $\varphi _{i}(n;k)\equiv 1$ then $\Phi (n)$
is the total number of particles born in the process within generations $%
0,1,...,n-1$; if $\varphi _{i}(n;k)=I\left\{ \sum_{j=1}^{m}\xi
_{ij}(n;k)\geq t\right\} $ for some positive integer $t$ (here and in what
follows $I\left\{ A\right\} $ means the indicator of the event $A$) and $%
\Phi (0)=0,$ then $\Phi (n)$ is the total number of particles of all types
in generations $0,1,...,n-1$ each of which had at least $t$ children, and so
on.

Letting $\lambda =0$ we arrive to the definition of the ordinary multitype
branching process in random environment (MBPRE) which we call the \textit{%
underlying} MBPRE for the initial MBPRCRE.

\begin{definition}
A $m$-type Galton-Watson process
\begin{equation*}
\mathbf{Z}(n)=\left( Z_{1}(n),\ldots ,Z_{m}(n)\right) ,\,n\in \mathbb{N}_{0}
\end{equation*}%
in a fixed (but selected at random) environment $\left\{ \mathbf{f}_{n}(%
\mathbf{s}),n\geq 0\right\} $ is a time-inhomogeneous Markov chain with the
state space
\begin{equation*}
\mathbb{N}_{0}^{m}:=\left\{ z=(z_{1},\ldots ,z_{m}),\,z_{i}\in \mathbb{N}%
_{0}\right\}
\end{equation*}%
defined as
\begin{equation}
\mathbf{Z}(0)=\mathbf{z},\;\mathbf{E}\left[ \mathbf{s}^{\mathbf{Z}(n+1)}|\;%
\mathbf{f}_{0},\ldots ,\mathbf{f}_{n};\mathbf{Z}(0),\ldots ,\mathbf{Z}(n)%
\right] =\left( \mathbf{f}_{n}\left( \mathbf{s}\right) \right) ^{\mathbf{Z}%
(n)}.  \label{defVP}
\end{equation}
\end{definition}

To simplify notation we write%
\begin{equation*}
\mathbf{E}_{\mathbf{f}}\left[ \mathbf{s}^{\mathbf{Z}(n+1)}\right] :=\mathbf{E%
}\left[ \mathbf{s}^{\mathbf{Z}(n+1)}|\;\mathbf{f}_{0},\ldots ,\mathbf{f}_{n}%
\right] .
\end{equation*}%
It follows from (\ref{defVP}) that $\mathbf{Z}(0)=\mathbf{z}$ and
\begin{equation*}
\mathbf{Z}(n+1):=\sum_{i=1}^{m}\sum_{k=1}^{Z_{i}(n)}\mathbf{\xi }_{i}(n;k)
\end{equation*}%
where $\mathbf{\xi }_{i}(n;k)\overset{d}{=}\mathbf{\xi }_{i}(n),k=1,2,\ldots
,Z_{i}(n),$ and, given the environment and $n=0,1,\ldots,$ and $%
i\in\{1,\ldots,m\}$ the mentioned random vectors are independent.

Let
\begin{equation}
A_{n}=\left( a_{ij}(n)\right) _{i,j=1}^{m}:=\left( \frac{\partial
F_{n}^{(i)}(\mathbf{s},\lambda )}{\partial s_{j}}\left\vert _{\mathbf{s}=%
\mathbf{1,}\lambda =0}\right. \right) _{i,j=1}^{m}=\left( \frac{\partial
f_{n}^{(i)}(\mathbf{s})}{\partial s_{j}}\left\vert _{\mathbf{s}=\mathbf{1}%
}\right. \right) _{i,j=1}^{m}  \label{meanMat}
\end{equation}%
be the mean matrix of the vector-valued p.g.f. $\mathbf{f}_{n}$ and%
\begin{equation}
\mathbf{C}_{n}:=\left( \mathbf{E}_{\mathbf{F}}\varphi _{1}(n),...,\mathbf{E}_{\mathbf{%
F}}\varphi _{1}(n)\right) ^{\prime }=\left( \frac{\partial F_{n}^{(1)}(%
\mathbf{s},\lambda )}{\partial \lambda }\left\vert _{\mathbf{s}=\mathbf{1,}%
\lambda =0}\right. ,...,\frac{\partial F_{n}^{(m)}(\mathbf{s},\lambda )}{%
\partial \lambda }\left\vert _{\mathbf{s}=\mathbf{1,}\lambda =0}\right.
\right) ^{\prime }.  \label{finprod}
\end{equation}%
By our assumptions the pairs $\left( A_{n},\mathbf{C}_{n}\right)
,\,n=0,1,... $ are iid: $\left( A_{n},\mathbf{C}_{n}\right)\overset{d}{=}( A,%
\mathbf{C}).$ Suppose that
\begin{equation}
\mathbf{E}\log ^{+}\left\Vert A\right\Vert <\infty .  \label{ExpCond}
\end{equation}%
It is known (see, for instance, \cite{Kin73}) that given condition \ (\ref%
{ExpCond}) the limit%
\begin{equation}
\lim_{n\rightarrow \infty }\frac{1}{n}\log \left\Vert A_{n-1}A_{n-2}\cdots
A_{0}\right\Vert =:\alpha  \label{LimMat0}
\end{equation}%
exists with probability 1 and, moreover,
\begin{equation}
\lim_{n\rightarrow \infty }\frac{1}{n}\mathbf{E}\log \left\Vert
A_{n-1}A_{n-2}\cdots A_{0}\right\Vert=\alpha .  \label{ExpecLim}
\end{equation}%
In what follows we call a MBPRE subcritical if $\alpha <0$ and
supercritical if $\alpha >0.$

Recall that the single-type BPRE with iid offspring p.g.f.'s were introduced
by Smith and Wilkinson in \cite{SW68} and, in a more general setting, in
\cite{AK71a} , \cite{AK71}, \cite{AN} and have been investigated by many
authors (see survey \cite{VZ} for a list of references up to 1985 and \cite%
{Koz}, \cite{VD2}, \cite{VD3},\cite{KOZ95}, \cite{AGKV}, \cite{AGKV2}, \cite%
{GDV}, \cite{GK00}, \cite{VD} and \cite{VD1} for some more recent results).
MBPRE were analyzed, in particular, in \cite{AK71a}, \cite{Key87b} and \cite%
{Tan81}.

Ordinary single-type Galton-Watson branching processes counted by random
characteristics where investigated by Sevastyanov \cite{Se72} ( for
integer-valued $\varphi (n;k)$) and by Grishechkin \cite{Gri90} (for the
general $\varphi (n;k)$). Grishechkin used the Galton-Watson and continuous
time Markov branching processes counted by random characteristics (with or
without immigration) to study queueing systems with processor sharing
discipline \cite{Gri88}.

\section{\protect\bigskip Limit theorems for MBPRCRE \label{SecLimT}}

Introduce the notation
\begin{equation*}
\Pi _{l,n}:=\prod_{i=l}^{n-1}A_{i},\qquad 1\leq l\leq n,
\end{equation*}%
with the agreement that $\Pi _{n,n}:=E$ is the unit $m\times m$ matrix. For
vectors $\mathbf{u}=(u_{1},\ldots ,u_{m}),\mathbf{v=}(v_{1},\ldots
,v_{m})^{\prime }\mathbf{\in }\mathbb{R}^{m}$ denote
\begin{equation*}
\left\langle \mathbf{u},\mathbf{v}\right\rangle :=\sum_{k=1}^{m}u_{i}v_{i}
\end{equation*}%
their inner product.

For a $m\times m$ \ matrix $A=\left( a_{ij}\right) _{i,j=1}^{m}$ and a $m-$%
dimensional vector $\mathbf{u=}(u_{1},\ldots ,u_{m})$  introduce
the norms
\begin{equation*}
\left\Vert A\right\Vert :=\sum_{i,j=1}^{m}\left\vert
a_{ij}\right\vert ,\quad \left\Vert \mathbf{u}\right\Vert
:=\sum_{i=1}^{m}\left\vert u_{i}\right\vert
\end{equation*}%
and%
\begin{equation*}
\left\Vert A\right\Vert _{2}:=\sqrt{\sum_{i,j=1}^{m}\left\vert
a_{ij}\right\vert ^{2}},\quad \left\Vert \mathbf{u}\right\Vert _{2}:=\sqrt{%
\sum_{i=1}^{m}\left\vert u_{i}\right\vert ^{2}}.
\end{equation*}

Now we formulate an important statement concerning properties of MBPRE.

Let%
\begin{equation*}
q_{i}\left( \mathbf{f}\right) :=\lim_{n\rightarrow \infty }\mathbf{P}_{%
\mathbf{f}}\left( \left\Vert \mathbf{Z}(n)\right\Vert =0|\mathbf{Z}(0)=%
\mathbf{e}_{i}\right)
\end{equation*}%
be the extinction probability of a MBPRE initiated at time $0$ by
a single individual of type $i$ and
\begin{equation*}
\mathbf{q}\left( \mathbf{f}\right) :=\left( q_{1}\left( \mathbf{f}\right)
,\ldots,q_{m}\left( \mathbf{f}\right) \right) .
\end{equation*}

\begin{theorem}
\label{Ttan}(\cite{Tan81}) If the mean matrices of a MBPRE meet condition (%
\ref{ExpCond}) and there exists a positive integer $L$ such that%
\begin{equation*}
\mathbf{P}\left( \min_{1\leq i,j\leq m}\left( A_{L-1}A_{L-2}\ldots
A_{0}\right) _{ij}>0\right) =1
\end{equation*}%
and $1\leq l\leq m$ such that%
\begin{equation*}
\mathbf{E}\left\vert \log \left( 1-\mathbf{P}_{\mathbf{f}}\left( \mathbf{Z}%
_{l}(L)=0|\mathbf{Z}(0)=\mathbf{e}_{i}\right) \right) \right\vert
<\infty ,
\end{equation*}%
then, for $\alpha $ specified by (\ref{ExpecLim})

1) $\alpha <0$ implies $\mathbf{P}\left( \mathbf{q}\left( \mathbf{f}\right) =%
\mathbf{1}\right) =1;$

2) $\alpha >0$ implies $\mathbf{P}\left( \mathbf{q}\left( \mathbf{f}\right) <%
\mathbf{1}\right) =1$ and%
\begin{equation}
\mathbf{P}_{\mathbf{f}}\left( \lim_{n\rightarrow \infty }n^{-1}\log
\left\Vert \mathbf{Z}(n)\right\Vert =\alpha |\mathbf{Z}(0)=\mathbf{e}%
_{i}\right) =1-q_{i}\left( \mathbf{f}\right)  \label{SU}
\end{equation}%
with probability 1 for $1\leq i\leq m.$
\end{theorem}

Let $\tau $ be the extinction moment of a MBPRE which starts by (may be
random) vector of the number of particles $\mathbf{Z}(0)$ with $\mathbf{E}%
\left\Vert \mathbf{Z}(0)\right\Vert <\infty $. Clearly,
\begin{eqnarray*}
\mathbf{P}\left( \tau >n\right) &=&\mathbf{P}\left( \left\Vert \mathbf{Z}%
(n)\right\Vert \geq 1\right) \leq \mathbf{E}\left\Vert \mathbf{Z}%
(n)\right\Vert =\mathbf{E}\left\Vert \mathbf{Z}(0)A_{0}A_{1}\cdots
A_{n-1}\right\Vert \\
&\leq &\mathbf{E}\left\Vert \mathbf{Z}(0)\right\Vert \mathbf{E}\left\Vert
A_{0}A_{1}\cdots A_{n-1}\right\Vert .
\end{eqnarray*}%
Hence we see that if $\alpha <0$ then for any $\alpha _{\ast }\in(0,-\alpha)
$ there exists a constant $K_{\ast }=K_{\ast }(\alpha _{\ast })\in (0,\infty
)$ such that for each $n=0,1,2,\ldots$
\begin{equation}
\mathbf{P}\left( \tau >n\right) \leq K_{\ast }e^{-\alpha _{\ast }n}.
\label{SurUp}
\end{equation}

Let
\begin{equation*}
\Phi :=\lim_{n\rightarrow \infty }\Phi (n)=\Phi (0)+\sum_{n=0}^{\infty
}\sum_{i=1}^{m}\sum_{k=1}^{Z_{i}(n)}\varphi _{i}(n;k)
\end{equation*}%
be the total size of the final product produced by the particles of the
MBPRCRE up to the extinction moment (if any). It is easy to see that if the
underlying\textbf{\ } MBPRE is supercritical, satisfies conditions of
Theorem \ref{Ttan} and, in addition, the final product of the MBPRCRE meets
the condition
\begin{equation}
\mathbf{P}\left( \min_{1\leq l\leq m}\mathbf{E}_{\mathbf{f}}\varphi
_{l}(n)>0\right) >0  \label{Nonzer}
\end{equation}%
then, by (\ref{SU}) and the law of large numbers, for each $i=1,\ldots ,m$
\begin{equation}
\mathbf{P}\left( \Phi =\infty |\mathbf{Z}(0)=\mathbf{e}_{i}\right) \geq
\mathbf{P}\left( \lim \inf_{n\rightarrow \infty }\left(
\sum_{j=1}^{m}\sum_{k=1}^{Z_{j}(n)}\varphi _{j}(n;k)\right) =\infty |\mathbf{%
Z}(0)=\mathbf{e}_{i}\right) >0.  \label{SupeInfin}
\end{equation}

Note, finally, that if $\Phi (0)=0$ and $\varphi _{j}(n;k)\equiv 1$ then
\begin{equation*}
\Phi :=\lim_{n\rightarrow \infty }\Phi (n)=\sum_{n=0}^{\infty }\left\Vert
\mathbf{Z}(n)\right\Vert
\end{equation*}%
is the total number of individuals ever existed in the MBPRCRE or, what is
the same, in the underlying MBPRE.

For a given $x\geq 0$ set%
\begin{equation}
s(x):=\lim_{n\rightarrow \infty }\left( \mathbf{E}\left\Vert A_{n-1}\cdots
A_{0}\right\Vert ^{x}\right) ^{1/n}=\lim_{n\rightarrow \infty }\left(
\mathbf{E}\left\Vert \Pi _{0,n}\right\Vert ^{x}\right) ^{1/n}  \label{defK}
\end{equation}%
and let%
\begin{equation}
s^{\prime }(0)=\lim_{n\rightarrow \infty }\frac{1}{n}\mathbf{E}\log
\left\Vert A_{n-1}\cdots A_{0}\right\Vert =\lim_{n\rightarrow \infty }\frac{1%
}{n}\mathbf{E}\log \left\Vert \Pi _{0,n}\right\Vert  \label{DefK0}
\end{equation}%
be the top Lyapunov exponent for this sequence of matrices.

Denote $\mathcal{D}:=\left\{ x>0:\mathbf{E}\left\Vert
A_{0}\right\Vert ^{x}<\infty \right\} $. It is known that the
limits in (\ref{defK}) and (\ref{DefK0})
exist and, moreover, $s(x)$ is a log-convex continuous function in $\mathcal{%
D}$ (see, for instance, \cite{Kin73}). Put
\begin{equation}
\kappa :=\inf \left\{ x>0:s(x)>1\right\}  \label{DefKappa}
\end{equation}%
and $\kappa =\infty $ if $s(x)\leq 1$ for all $x>0$. Observe that $s(0)=1$
and, therefore, $\kappa =0$ if $s^{\prime }(0)>0$ and $\kappa \in (0,\infty
] $ if $s^{\prime }(0)<0$.

In the last case (which will be our main concern) the series $%
\sum_{n=0}^{\infty }\mathbf{E}\left\Vert \Pi _{0,n}\right\Vert ^{x}$
converges if $0<x<\kappa $ and diverges if $x>\kappa $.

Introduce the set
\begin{equation*}
U_{+}=\left\{ \mathbf{u}=(u_{1},\ldots ,u_{m})\in \mathbb{R}^{m}:u_{i}\geq
0,1\leq i\leq m,\,\left\Vert \mathbf{u}\right\Vert _{2}=1\right\}
\end{equation*}%
and associate with the tuple $(A_{n},\mathbf{C}_{n}),n=0,1,2,\ldots $ of iid
pairs the series
\begin{equation*}
\Xi _{l}:=\sum_{k=l}^{\infty }A_{l}A_{l+1}\ldots A_{k-1}\mathbf{C}%
_{k}=\sum_{k=l}^{\infty }\Pi _{l,k}\mathbf{C}_{k},\ l=0,1,\ldots ;\,\quad
\Xi =:\Xi _{0}.
\end{equation*}

Our main results are established under the following hypothesis.

\textbf{Condition T.} There exist positive constants $\kappa $ and $K_{0}$
and a continuous strictly positive function $l(\mathbf{u})$ on $U_{+}$ such
that for all $\mathbf{u\in }U_{+}$
\begin{equation*}
\lim_{y\rightarrow \infty }y^{\kappa }\mathbf{P}\left( \left\langle \mathbf{%
u,}\Xi \right\rangle >y\right) =K_{0}l(\mathbf{u}).
\end{equation*}

In Section \ref{SecAux} we list sufficient conditions on the
distributions of the pairs $(A_{n},\mathbf{C}_{n})$ which provide
the validity of Condition $T$. These conditions are extracted from
paper \cite{Kes74} where the behavior of the tail distribution of
sums and products of random matrices were investigated.

The following theorem is the main result of the article.

\begin{theorem}
\label{TailMain}Let a MBPRCRE satisfy the following hypotheses:

1) the underlying MBPRE \ is subcritical and meet conditions of Theorem \ref%
{Ttan};

2) for $\kappa $ specified by (\ref{DefKappa}) the following assumptions
fulfill:

if $\kappa >1$ then
\begin{equation}
\max_{1\leq i\leq m}\mathbf{E}\left\vert \sum_{j=1}^{m}\left( \xi
_{ij}-a_{ij}\right) \right\vert ^{\kappa }<\infty \text{ and \ \ }\mathbf{E}%
\left\vert \sum_{i=1}^{m}\left( \varphi _{i}(n)-\mathbf{E}_{\mathbf{F}%
}\varphi _{i}(n)\right) \right\vert ^{\kappa }<\infty ,  \label{Ckappabig}
\end{equation}

if $\kappa \leq 1$ then%
\begin{equation}
\max_{1\leq i\leq m}\mathbf{E}\left( \sum_{j=1}^{m}Var_{\mathbf{F}}\xi
_{ij}\right) ^{\kappa }<\infty \text{ and \ }\mathbf{E}\left(
\sum_{i=1}^{m}Var_{\mathbf{F}}\varphi _{i}(n)\right) ^{\kappa }<\infty ;
\label{Ckappasmall}
\end{equation}

4) there exists $\delta >0$ such that $0<\mathbf{E}\varphi _{i}^{\kappa
+\delta }(n)<\infty ,\quad i=1,\ldots ,m. $

If, in addition, the mean matrix (\ref{meanMat}) and the vector (\ref%
{finprod}) are such that Condition $T$ is valid, then, as $y\rightarrow
\infty $%
\begin{equation*}
\mathbf{P}\left( \Phi >y\right) \sim Cy^{-\kappa },\qquad C\in (0,\infty ).
\end{equation*}
\end{theorem}

An evident corollary of Theorem \ref{TailMain} is the following statement.

\begin{theorem}
\label{Tmoment}If the conditions of Theorem \ref{TailMain} are valid then
\begin{equation*}
\mathbf{E}\Phi ^{x}<\infty
\end{equation*}%
if and only if $x<\kappa .$
\end{theorem}

\section{Auxiliary results\label{SecAux}}

The proof of Theorem \ref{TailMain} is heavily based on Condition $T$ whose
validity is not easy to check. We list here a set of assumptions given in
\cite{Kes74} which imply Condition T.

Let $\Lambda \left( A\right) $ be the spectral radius of the matrix $A$. The
following statement is a refinement of a Kesten theorem from \cite{Kes74}.

\begin{theorem}
\label{Tkest}(see \cite{DMB99}) Let $\left\{ A_{n},n\geq 0\right\} $ be a
sequence of iid matrices generated by a measure $\mathbb{P}_{A}$ with
support concentrated on nonnegative matrices and $A=\left( a_{ij}\right)
_{i,j=1}^{m}\overset{d}{=}A_{n}.$ Assume that the following conditions are
valid:

1) there exists $\varepsilon >0$ such that $E\left\Vert A\right\Vert
^{\varepsilon }<\infty $;

2) $A$ has no zero rows a.s.;

3) the group generated by
\begin{equation*}
\left\{ \log \Lambda (a_{n}\cdots a_{0}):a_{n}\cdots a_{0}>0\text{ for some }%
n\text{ and }a_{i}\in \text{supp}(\mathbb{P}_{A})\right\}
\end{equation*}%
is dense in $\mathbb{R};$

4) there exists $\kappa _{0}>0$ for which
\begin{equation*}
\mathbf{E}\left[ \min_{1\leq i\leq m}\left( \sum_{j=1}^{m}a_{ij}\right)
^{\kappa _{0}}\right] \geq m^{\kappa _{0}/2}
\end{equation*}%
and%
\begin{equation*}
\mathbf{E}\left\Vert A\right\Vert ^{\kappa _{0}}\log ^{+}\left\Vert
A\right\Vert <\infty .
\end{equation*}

Then there exists a $\kappa \in (0,\kappa _{0}]$ such that%
\begin{equation*}
s^{\prime }(\kappa )=\lim_{n\rightarrow \infty }\frac{1}{n}\log \mathbf{E}%
\left\Vert A_{n-1}\cdots A_{0}\right\Vert ^{\kappa }=0.
\end{equation*}

If, in addition, the tuple of $m$-dimensional vectors
$\left\{\mathbf{C}_{n},n\geq 0\right\}$
is such that the pairs $(A_{n},\mathbf{C}_{n})$, $n=0,1,\ldots $ are iid: $%
\left( A_{n},\mathbf{C}_{n}\right) \overset{d}{=}\left( A,\mathbf{C}\right) $
and such that
\begin{equation*}
\mathbf{P}\left( \mathbf{C}=\mathbf{0}\right) <1,\quad \mathbf{P}\left(
\mathbf{C}\geq \mathbf{0}\right) =1,\quad \mathbf{E}\left\Vert \mathbf{C}%
\right\Vert ^{\kappa }<\infty ,
\end{equation*}%
then there exist a constant $K_{0}\in \left( 0,\infty \right) $ and a
continuous strictly positive function $l(\mathbf{u})$ on $U_{+}$ such that%
\begin{equation*}
\lim_{y\rightarrow \infty }y^{\kappa }\mathbf{P}\left( \left\langle \mathbf{%
u,}\Xi \right\rangle >y\right) =K_{0}l(\mathbf{u}),\ \mathbf{u\in }U_{+}.
\end{equation*}
\end{theorem}

The next lemma will be of importance for subsequent arguments.

\begin{lemma}
\label{LGut}(\cite{Gut}, Theorem 1.5.1) If $X_{i},i=1,2,\ldots $ is a
sequence of iid random variables such that $\mathbf{E}\left\vert
X_{i}\right\vert ^{p}<\infty $ and $\mathbf{E}X_{i}=0$ if $p\geq 1,$ and $N$
is a stopping time for the sequence $S_{n}=X_{1}+\ldots +X_{n},$ then there
exists a constant $R_{p}\in (0,\infty )$ such that
\begin{equation*}
\mathbf{E}\left\vert S_{N}\right\vert ^{p}\leq R_{p}\mathbf{E}\left\vert
X_{i}\right\vert ^{p}\mathbf{E}N^{p/2\vee 1}.
\end{equation*}
\end{lemma}

\section{Properties of the underlying MBPRE\label{NonRand}}

In this section we assume that the conditions of Theorem
\ref{TailMain} are valid. This means, in particular, that we deal
with subcritical MBPRE. Let us agree to denote by $K,K_{x},x\in
(0,\infty )$ positive constants which may be different from
formula to formula.

First
we evaluate the expectation of the random variable $\left\Vert \mathbf{Z}%
(n)\right\Vert ^{x},\,0<x<\kappa ,$ from above.
\begin{lemma}
\label{LGeom}If $\mathbf{Z}(0)=\mathbf{z}$ and $\kappa >1$ then for each $%
x\in \lbrack 1,\kappa )$ there exist $\rho _{x}\in \left( 0,1\right) $ and $%
K_{x}<\infty $ such that%
\begin{equation}
\mathbf{E}\left\Vert \mathbf{Z}(n)\right\Vert ^{x}\leq K_{x}\rho
_{x}^{n}\left\Vert \mathbf{z}\right\Vert ^{x}  \label{geom11}
\end{equation}%
for all $n=1,2,\ldots .$
\end{lemma}

\textbf{Proof}. Clearly, for any nonrandom vector $\mathbf{b}\in \mathbb{R}%
^{m}$
\begin{eqnarray*}
\mathbf{E}\left\langle \mathbf{Z}(n),\mathbf{b}\right\rangle &=&\mathbf{E}%
\left\langle \sum_{i=1}^{m}\sum_{k=1}^{Z_{i}(n-1)}\mathbf{\xi }_{i}(n-1;k),%
\mathbf{b}\right\rangle \\
&=&\mathbf{E}\left\langle \mathbf{Z}(n-1),A_{n-1}\mathbf{b}\right\rangle =%
\mathbf{E}\left\langle \mathbf{Z}(n-j),\Pi _{n-j,n}\mathbf{b}\right\rangle \\
&=&\ldots =\mathbf{E}\left\langle \mathbf{z,}\Pi _{0,n}\mathbf{b}%
\right\rangle ,\qquad 0\leq j\leq n.
\end{eqnarray*}%
Thus, if $\left\Vert \mathbf{b}\right\Vert \leq $ $K<\infty $ then for any $%
\delta >0$ with $s(1)\left( 1+\delta \right) <1$ there exists a constant $%
K_{1}\in \left( 0,\infty \right) $ such that for all $n=1,2,\ldots $ and $%
0\leq j\leq n$%
\begin{equation}
\left\vert \mathbf{E}\left\langle \mathbf{Z}(n-j),\Pi _{n-j,n}\mathbf{b}%
\right\rangle \right\vert \leq \left\Vert \mathbf{z}\right\Vert \mathbf{E}%
\left\Vert \Pi _{0,n}\right\Vert \mathbf{E}\left\Vert \mathbf{b}\right\Vert
\leq K_{1}\left( s(1)(1+\delta )\right) ^{n}\left\Vert \mathbf{z}\right\Vert
.  \label{geom12}
\end{equation}

Now we apply arguments similar to those used in \cite{DPZ96}. It is easy to
check that for any $y,w\geq 0$ and any $\varepsilon \in (0,1)$%
\begin{equation*}
y^{x}\leq (1+\varepsilon )w^{x}+c_{x,\varepsilon }\left\vert y-w\right\vert
^{x},
\end{equation*}%
where $c_{x,\varepsilon }:=\left( 1-\left( 1+\varepsilon \right)
^{-1/x}\right) ^{-x}$. Hence we have
\begin{eqnarray}
\mathbf{E}\left\Vert \mathbf{Z}(n)\right\Vert ^{x} &=&\mathbf{E}\left\Vert
\mathbf{Z}(n-1)A_{n-1}+\mathbf{Z}(n)-\mathbf{Z}(n-1)A_{n-1}\right\Vert ^{x}
\notag \\
&\leq &(1+\varepsilon )\mathbf{E}\left\Vert \mathbf{Z}(n-1)A_{n-1}\right%
\Vert ^{x}+c_{x,\varepsilon }\mathbf{E}\left\Vert \mathbf{Z}(n)-\mathbf{Z}%
(n-1)A_{n-1}\right\Vert ^{x}.  \notag \\
&&  \label{SSta}
\end{eqnarray}%
Recalling the definition $a_{ij}(n)=\mathbf{E}_{\mathbf{F}}\xi _{ij}(n)$, set%
\begin{equation}
\beta _{i}(n):=\sum_{j=1}^{m}(\mathbf{\xi }_{ij}(n)-a_{ij}(n))
\label{DefetaBeta}
\end{equation}%
and let%
\begin{equation*}
M_{x}(n;i):=\mathbf{E}_{\mathbf{F}}\left\vert \beta _{i}(n)\right\vert
^{x},\ \ M_{x}:=\max_{1\leq i\leq m}\mathbf{E}\left\vert \beta
_{i}(n)\right\vert ^{x}.
\end{equation*}%
By Lemma \ref{LGut} with $p=x>1$ we conclude%
\begin{align}
& \mathbf{E}_{\mathbf{F}}\left\Vert \mathbf{Z}(n)-\mathbf{Z}%
(n-1)A_{n-1}\right\Vert ^{x}  \notag \\
=& \mathbf{E}_{\mathbf{F}}\left\vert
\sum_{i=1}^{m}\sum_{k=1}^{Z_{i}(n-1)}\sum_{j=1}^{m}\left[ \mathbf{\xi }%
_{ij}(n-1;k)-a_{ij}(n-1)\right] \right\vert ^{x}  \notag \\
& \leq m^{x}\sum_{i=1}^{m}\mathbf{E}_{\mathbf{F}}\left\vert
\sum_{k=1}^{Z_{i}(n-1)}\sum_{j=1}^{m}\left[ \mathbf{\xi }%
_{ij}(n-1;k)-a_{ij}(n-1)\right] \right\vert ^{x}  \notag \\
& \leq R_{x}m^{x}\sum_{i=1}^{m}M_{x}(n-1,i)\mathbf{E}_{\mathbf{F}%
}Z_{i}^{x/2\vee 1}(n-1).  \label{Mom22}
\end{align}%
These estimates and (\ref{SSta}) give
\begin{eqnarray}
\mathbf{E}\left\Vert \mathbf{Z}(n)\right\Vert ^{x} &\leq &(1+\varepsilon )%
\mathbf{E}\left\Vert \mathbf{Z}(n-1)A_{n-1}\right\Vert ^{x}  \notag \\
&&+c_{x,\varepsilon }R_{x}m^{x}M_{x}\sum_{i=1}^{m}\mathbf{E}Z_{i}^{x/2\vee
1}(n-1).  \label{above21}
\end{eqnarray}

Now we are ready to demonstrate (\ref{geom11}) . First we establish (\ref%
{geom11}) for all integer $x\in \lbrack 1,\kappa )$. For $x=1$ we have
proved (\ref{geom11}) by (\ref{geom12}) with $\mathbf{b}=\mathbf{1}$. Now we
use induction for $x\geq 2$. Observing that $x/2\vee 1\leq x-1$ in this
case, we see by (\ref{above21}) and the estimate
\begin{equation*}
\sum_{i=1}^{m}\mathbf{E}Z_{i}^{x/2\vee 1}(n-1)\leq \mathbf{E}%
\sum_{i=1}^{m}Z_{i}^{x-1}(n-1)\leq \mathbf{E}\left\Vert \mathbf{Z}%
(n-1)\right\Vert ^{x-1}
\end{equation*}%
that%
\begin{eqnarray}
\mathbf{E}\left\Vert \mathbf{Z}(n)\right\Vert ^{x} &\leq &(1+\varepsilon )%
\mathbf{E}\left\Vert \mathbf{Z}(n-1)A_{n-1}\right\Vert ^{x}  \notag \\
&&+c_{x,\varepsilon }R_{x}m^{x}M_{x}\mathbf{E}\left\Vert \mathbf{Z}%
(n-1)\right\Vert ^{x-1} \leq (1+\varepsilon )^{2}\mathbf{E}\left\Vert
\mathbf{Z}(n-2)A_{n-2}A_{n-1}\right\Vert ^{x}  \notag \\
&&+(1+\varepsilon )c_{x,\varepsilon }R_{x}m^{x}M_{x}\mathbf{E}\left\Vert
A_{n-1}\right\Vert ^{x}\mathbf{E}\left\Vert \mathbf{Z}(n-2)\right\Vert ^{x-1}
\notag \\
&&+c_{x,\varepsilon }R_{x}m^{x}M_{x}\mathbf{E}\left\Vert \mathbf{Z}%
(n-1)\right\Vert ^{x-1} \leq \ldots \leq (1+\varepsilon )^{n}\mathbf{E}%
\left\Vert \mathbf{z}\Pi _{0,n}\right\Vert ^{x}  \notag \\
&&+c_{x,\varepsilon }R_{x}m^{x}M_{x}\sum_{j=0}^{n-1}(1+\varepsilon )^{j}%
\mathbf{E}\left\Vert \Pi _{n-j,n}\right\Vert ^{x}\mathbf{E}\left\Vert
\mathbf{Z}(n-j-1)\right\Vert ^{x-1}.  \notag \\
&&  \label{MMM11}
\end{eqnarray}

Since $x\in [1,\kappa)$, for any $\delta >0$ there exists a constant $L_{x}$
such that
\begin{equation*}
\mathbf{E}\left\Vert \Pi _{0,n}\right\Vert ^{x}\leq L_{x}\left(
s(x)(1+\delta )\right) ^{n}
\end{equation*}%
for all $n=0,1,2,\ldots .$ By induction hypothesis there exist constants $%
K_{x-1}$ and $\rho _{x-1}\in (s(x),1)$ such that
\begin{equation*}
\mathbf{E}\left\Vert \mathbf{Z}(n-j)\right\Vert ^{x-1}\leq K_{x-1}\rho
_{x-1}^{n-j}\left\Vert \mathbf{z}\right\Vert ^{x-1}\leq K_{x-1}\rho
_{x-1}^{n-j}\left\Vert \mathbf{z}\right\Vert ^{x}
\end{equation*}%
for all $j=0,1,\ldots ,n$. Thus,%
\begin{eqnarray}
&&\mathbf{E}\left\Vert \mathbf{Z}(n)\right\Vert ^{x}\leq (1+\varepsilon
)^{n}L_{x}\left( s(x)(1+\delta )\right) ^{n}\left\Vert \mathbf{z}\right\Vert
^{x}  \notag \\
&&+c_{x,\varepsilon }R_{x}m^{x}M_{x}L_{x}K_{x-1}\left\Vert \mathbf{z}%
\right\Vert ^{x}\sum_{j=0}^{n-1}(1+\varepsilon )^{j}\left( s(x)(1+\delta
)\right) ^{j}\rho _{x-1}^{n-j-1}  \notag \\
&&\quad =L_{x}s^{n}(x)(1+\varepsilon )^{n}\left( 1+\delta \right)
^{n}\left\Vert \mathbf{z}\right\Vert ^{x}  \notag \\
&&\qquad +Ks^{n}(x)(1+\varepsilon )^{n}\left( 1+\delta \right)
^{n}\left\Vert \mathbf{z}\right\Vert ^{x}\sum_{j=1}^{n}\frac{\rho
_{x-1}^{j-1}}{s^{j}(x)(1+\varepsilon )^{j}\left( 1+\delta \right) ^{j}}.
\label{ImmigMathem1}
\end{eqnarray}%
Now selecting $\delta $ and $\varepsilon $ in such a way that%
\begin{equation*}
s(x)(1+\varepsilon )\left( 1+\delta \right) \in \left( \rho _{x-1},1\right)
\end{equation*}%
(the last is always possible since $\rho _{x-1}<s(x)<1$) we get (\ref{geom11}%
) .

To treat the case of noninteger $x\in \lbrack 1,\kappa )$ observe that $%
x^{\ast }=\left[ \kappa \right] \in \lbrack 1,\kappa )$ is an integer for
which (\ref{geom11}) \ is valid. Thus, it remains to demonstrate (\ref%
{geom11}) for $x=x^{\ast }+\gamma <\kappa ,$ where $\gamma \in (0,1)$. Since
$x^{\ast }\geq x/2\vee 1$ one can use the same arguments as earlier with $%
x^{\ast }$ for $x-1$.

The lemma is proved.

Let $r$ be an integer and
\begin{equation*}
\zeta =\zeta (r):=\min \left\{ n\geq 0:\,\left\Vert \mathbf{Z}(n)\right\Vert
>r\right\}
\end{equation*}%
with the natural agreement that $\zeta =\infty $ if $\max_{n}\left\Vert
\mathbf{Z}(n)\right\Vert \leq r$.

\begin{lemma}
\label{Lke1} Under the conditions of Theorem \ref{TailMain} for any fixed $%
r\geq 1$%
\begin{equation*}
\mathbf{E}\left[ \left\Vert \mathbf{Z}(\zeta (r))\right\Vert ^{\kappa
}I\left\{ \zeta (r)<\infty \right\} \right] <\infty .
\end{equation*}
\end{lemma}

\textbf{Proof}. Similar to (\ref{above21}) we have for $x>1$
\begin{eqnarray}
&&\mathbf{E}_{\mathbf{F}}\left[ \left\Vert \mathbf{Z}(n)\right\Vert ^{x}|%
\mathbf{Z}(0),\ldots ,\mathbf{Z}(n-1)\right] \leq (1+\varepsilon )\left\Vert
\mathbf{Z}(n-1)A_{n-1}\right\Vert ^{x}  \notag \\
&&\qquad +c_{x,\varepsilon
}R_{x}m^{x}\sum_{i=1}^{m}M_{x}(n-1,i)Z_{i}^{x/2\vee 1}(n-1)=:\Psi _{x}(n-1),
\label{PP1}
\end{eqnarray}%
while Jensen's inequality yields for $x\leq 1$
\begin{eqnarray}
\mathbf{E}_{\mathbf{F}}\left[ \left\Vert \mathbf{Z}(n)\right\Vert ^{x}|%
\mathbf{Z}(0),\ldots ,\mathbf{Z}(n-1)\right] &\leq &\left( \mathbf{E}_{%
\mathbf{F}}\left[ \left\Vert \mathbf{Z}(n)\right\Vert |\mathbf{Z}(0),\ldots ,%
\mathbf{Z}(n-1)\right] \right) ^{x}  \notag \\
&&  \notag \\
&=&\left\Vert \mathbf{Z}(n-1)A_{n-1}\right\Vert ^{x}\leq \Psi _{x}(n-1).
\label{PP3}
\end{eqnarray}%
Clearly,
\begin{eqnarray*}
\Psi _{\kappa}(n-1)I\left\{ n<\zeta \right\} &\leq
&Q_{n-1}(r):=(1+\varepsilon )r^{\kappa }\left\Vert A_{n-1}\right\Vert
^{\kappa } \\
&&+c_{x,\varepsilon }R_{x}m^{\kappa }r^{\kappa /2\vee
1}\sum_{i=1}^{m}M_{\kappa }(\zeta -1,i).
\end{eqnarray*}%
Using these estimates we have on the event $\left\{ \zeta <\infty \right\} :$
\begin{eqnarray}
\left\Vert \mathbf{Z}(\zeta )\right\Vert ^{\kappa }=\Psi _{\kappa }(\zeta -1)%
\frac{\left\Vert \mathbf{Z}(\zeta )\right\Vert ^{\kappa }}{\Psi _{\kappa
}(\zeta -1)} &\leq &Q_{\zeta -1}(r)\frac{\left\Vert \mathbf{Z}(\zeta
)\right\Vert ^{\kappa }}{\Psi _{\kappa }(\zeta -1)}  \notag \\
&\leq &\sum_{\zeta \leq n<\tau }Q_{n}(r)\frac{\left\Vert \mathbf{Z}%
(n)\right\Vert ^{\kappa }}{\Psi _{\kappa }(n-1)}.  \label{ko1}
\end{eqnarray}%
By (\ref{PP1}) and (\ref{PP3}) we now see that%
\begin{eqnarray*}
&&\mathbf{E}\left[ \left\Vert Z(\zeta )\right\Vert ^{\kappa }I\left\{ \zeta
<\infty \right\} \right] \leq \sum_{n\geq 1}\mathbf{E}\left[ Q_{n}(r)\frac{%
\left\Vert \mathbf{Z}(n)\right\Vert ^{\kappa }}{\Psi _{\kappa }(n-1)}%
I\left\{ \tau \geq n\right\} \right] \\
&&\quad \leq \sum_{n\geq 1}\mathbf{E}\left[ Q_{n}(r)I\left\{ \tau \geq
n\right\} \right] =\sum_{n\geq 1}\mathbf{P}\left( \tau \geq n\right) \mathbf{%
E}\left[ Q_{n}(r)\right] \\
&&\qquad =\left[ (1+\varepsilon )r^{\kappa }\mathbf{E}\left\Vert
A\right\Vert ^{\kappa }+c_{x,\varepsilon }m^{\kappa +1}r^{\kappa /2\vee
1}M_{\kappa }\right] \sum_{n\geq 1}\mathbf{P}\left( \tau \geq n\right)
<\infty ,
\end{eqnarray*}%
since for $\kappa _{1}:=\min \left( \kappa /2,1\right) $%
\begin{eqnarray*}
\sum_{n\geq 1}\mathbf{P}\left( \tau \geq n\right) &=&\sum_{n\geq 1}\mathbf{P}%
\left( \left\Vert \mathbf{Z}(n-1)\right\Vert \geq 1\right) =\sum_{n\geq 1}%
\mathbf{E}\left[ \mathbf{P}_{\mathbf{F}}\left( \left\Vert \mathbf{Z}%
(n-1)\right\Vert \geq 1\right) \right] \\
&\leq &\sum_{n\geq 1}\mathbf{E}\left[ \mathbf{E}_{\mathbf{F}}\left\Vert
\mathbf{Z}(n-1)\right\Vert ^{\kappa _{1}}\right] \leq \sum_{n\geq 1}\mathbf{E%
}\left\Vert \mathbf{E}_{\mathbf{F}}\mathbf{Z}(n-1)\right\Vert ^{\kappa _{1}}
\\
&\leq &\sum_{n\geq 1}\left\Vert z\right\Vert ^{\kappa _{1}}\mathbf{E}%
\left\Vert A_{0}\cdots A_{n-1}\right\Vert ^{\kappa _{1}}<\infty .
\end{eqnarray*}

The lemma is proved.

Let $\mathcal{B}_{n}$ ($n=1,2,\ldots $) be the $\sigma $-algebra generated
by the tuple%
\begin{equation*}
\mathbf{F}_{0}(\mathbf{s},\lambda ),\mathbf{F}_{1}(\mathbf{s},\lambda
),\ldots ,\mathbf{F}_{n-1}(\mathbf{s},\lambda ),\mathbf{Z}(0),\ldots ,%
\mathbf{Z}(n)
\end{equation*}%
and let
\begin{equation*}
\mathbf{C}_{n}:=\left( \mathbf{E}_{\mathbf{F}}\varphi _{1}(n),...,\mathbf{E}%
_{\mathbf{F}}\varphi _{m}(n)\right) ^{\prime },\,n=0,1,2,\ldots
\end{equation*}

Recall that the pairs $\left( A_{n},\mathbf{C}_{n}\right) ,n=0,1,\ldots ,$
are iid according to the definition of our MBPRCRE and, in particular, $%
\left( A_{n},\mathbf{C}_{n}\right) $ is independent on $\mathcal{B}_{n}$ .
Set
\begin{equation}
S(\zeta ):=\sum_{n=\zeta }^{\infty }\left\langle \mathbf{Z}(n),\mathbf{C}%
_{n}\right\rangle .\qquad  \label{Defeta}
\end{equation}

The next lemma shows that for large $r$ the random variable $S(\zeta
)=S(\zeta (r))$ is, in a sense, close to the conditional expectation $%
\mathbf{E}\left[ S(\zeta )|\mathcal{B}_{\zeta }\right] =\left\langle \mathbf{%
Z}(\zeta ),\Xi _{\zeta }\right\rangle .$

\begin{lemma}
\label{Lke2}Under the conditions of Theorem \ref{TailMain} for any $%
\varepsilon >0$ there exists $r=r(\varepsilon )$ such that for all $y\geq
y_{0}$%
\begin{equation}
\mathbf{P}\left( \left\vert S(\zeta )-\left\langle \mathbf{Z}(\zeta ),\Xi
_{\zeta }\right\rangle \right\vert >\varepsilon y;\zeta <\infty \right) \leq
\frac{\varepsilon }{y^{\kappa }}\mathbf{E}\left[ \left\Vert \mathbf{Z}(\zeta
)\right\Vert ^{\kappa }I\left\{ \zeta <\infty \right\} \right] .
\label{Diff1}
\end{equation}
\end{lemma}

\textbf{Proof}. Evidently, for $n\geq \zeta +1$%
\begin{eqnarray*}
\left\langle \mathbf{Z}(n)-\mathbf{Z}(\zeta )\Pi _{\zeta ,n},\mathbf{C}%
_{n}\right\rangle &=&\sum_{l=\zeta +1}^{n}\left\langle \mathbf{Z}(l)\Pi
_{l,n}-\mathbf{Z}(l-1)\Pi _{l-1,n},\mathbf{C}_{n}\right\rangle \\
&=&\sum_{l=\zeta +1}^{n}\left\langle \mathbf{Z}(l)-\mathbf{Z}%
(l-1)A_{l-1},\Pi _{l,n}\mathbf{C}_{n}\right\rangle
\end{eqnarray*}%
which implies%
\begin{eqnarray}
\left\vert S(\zeta )-\left\langle \mathbf{Z}(\zeta ),\Xi _{\zeta
}\right\rangle \right\vert &=&\left\vert \sum_{n=\zeta +1}^{\infty
}\left\langle \mathbf{Z}(n)-\mathbf{Z}(\zeta )\Pi _{\zeta ,n},\mathbf{C}%
_{n}\right\rangle \right\vert  \notag \\
&\leq &\sum_{n=\zeta +1}^{\infty }\left\vert \left\langle \mathbf{Z}(n)-%
\mathbf{Z}(\zeta )\Pi _{\zeta ,n},\mathbf{C}_{n}\right\rangle \right\vert
\notag \\
&=&\sum_{n=\zeta +1}^{\infty }\left\vert \sum_{l=\zeta +1}^{n}\left\langle
\mathbf{Z}(l)-\mathbf{Z}(l-1)A_{l-1},\Pi _{l,n}\mathbf{C}_{n}\right\rangle
\right\vert  \notag \\
&\leq &\sum_{l=\zeta +1}^{\infty }\sum_{n=l}^{\infty }\left\vert
\left\langle \mathbf{Z}(l)-\mathbf{Z}(l-1)A_{l-1},\Pi _{l,n}\mathbf{C}%
_{n}\right\rangle \right\vert  \notag \\
&\leq &\sum_{l=\zeta +1}^{\infty }\left\Vert \mathbf{Z}(l)-\mathbf{Z}%
(l-1)A_{l-1}\right\Vert \left\Vert \sum_{n=l}^{\infty }\Pi _{l,n}\mathbf{C}%
_{n}\right\Vert  \notag \\
&=&\sum_{l=\zeta +1}^{\infty }\left\Vert \mathbf{Z}(l)-\mathbf{Z}%
(l-1)A_{l-1}\right\Vert \left\Vert \Xi _{l}\right\Vert .  \label{seco2}
\end{eqnarray}

Hence, on account of $\sum_{j=1}^{\infty }j^{-2}=\pi ^{2}/6\leq 2,$ we may
apply on the event $\{ \zeta <\infty\}$ the arguments used in \cite{KKS75},
Lemma 3, to conclude that
\begin{eqnarray}
&&\mathbf{P}_{\mathbf{F}}\left( \left\vert S(\zeta )-\left\langle \mathbf{Z}%
(\zeta ),\Xi _{\zeta }\right\rangle \right\vert \geq \varepsilon
y\,\left\vert \,\mathcal{B}_{\zeta }\right. \right)  \notag \\
&\leq &\mathbf{P}_{\mathbf{F}}\left( \sum_{l=\zeta +1}^{\infty }\left\Vert
\mathbf{Z}(l)-\mathbf{Z}(l-1)A_{l-1}\right\Vert \left\Vert \Xi
_{l}\right\Vert \geq 6\pi ^{-2}\varepsilon y\sum_{l=\zeta +1}^{\infty }\frac{%
1}{\left( l-\zeta \right) ^{2}}\left\vert \,\mathcal{B}_{\zeta }\right.
\right)  \notag \\
&&\quad \leq \sum_{l=\zeta +1}^{\infty }\mathbf{P}_{\mathbf{F}}\left(
\left\Vert \mathbf{Z}(l)-\mathbf{Z}(l-1)A_{l-1}\right\Vert \left\Vert \Xi
_{l}\right\Vert \geq \frac{\varepsilon y}{2\left( l-\zeta \right) ^{2}}%
\,\left\vert \mathcal{B}_{\zeta }\right. \right) .  \label{BBB1}
\end{eqnarray}%
Since $\mathbf{Z}(l)-\mathbf{Z}(l-1)A_{l-1}$ and $\Xi _{l}$ are independent
random vectors on the event $\zeta \leq l<\infty $, we get%
\begin{eqnarray*}
&&\mathbf{P}_{\mathbf{F}}\left( \left\Vert \mathbf{Z}(l)-\mathbf{Z}%
(l-1)A_{l-1}\right\Vert \left\Vert \Xi _{l}\right\Vert \geq \frac{%
\varepsilon y}{2\left( l-\zeta \right) ^{2}}\,\left\vert \,\mathcal{B}%
_{\zeta }\right. \right) \\
&&\quad =\int_0^\infty  \mathbf{P}_{\mathbf{F}}\left( \left\Vert \mathbf{Z}(l)-\mathbf{%
Z}(l-1)A_{l-1}\right\Vert \in dt\,\left\vert \,\mathcal{B}_{\zeta
}\right.
\right) \mathbf{P}\left( \left\Vert \Xi \right\Vert \geq \frac{\varepsilon y%
}{2t\left( l-\zeta \right) ^{2}}\right) .
\end{eqnarray*}%
According to Condition $T$ there exists a constant $K\in \left( 0,\infty
\right) $ such that for all $l>\zeta $
\begin{eqnarray}
&&\mathbf{P}_{\mathbf{F}}\left( \left\Vert \mathbf{Z}(l)-\mathbf{Z}%
(l-1)A_{l-1}\right\Vert \left\Vert \Xi _{l}\right\Vert \geq \frac{%
\varepsilon y}{2\left( l-\zeta \right) ^{2}}\,\left\vert \mathcal{B}_{\zeta
}\right. \right)  \notag \\
&&\quad \leq \int_0^\infty  \mathbf{P}_{\mathbf{F}}\left( \left\Vert \mathbf{Z}(l)-%
\mathbf{Z}(l-1)A_{l-1}\right\Vert \in dt\,\left\vert
\mathcal{B}_{\zeta
}\right. \right) \frac{Kt^{\kappa }}{\varepsilon ^{\kappa }y^{\kappa }}%
\left( l-\zeta \right) ^{2\kappa }  \notag \\
&&\qquad \leq \frac{K}{\varepsilon ^{\kappa }y^{\kappa }}\left( l-\zeta
\right) ^{2\kappa }\mathbf{E}_{\mathbf{F}}\left[ \left\Vert \mathbf{Z}(l)-%
\mathbf{Z}(l-1)A_{l-1}\right\Vert ^{\kappa }\left\vert \,\mathcal{B}_{\zeta
}\right. \right] .  \label{BB}
\end{eqnarray}%
Now we consider the cases $\kappa \leq 1$ and $\kappa >1$ separately.

For the first case we use the estimate%
\begin{equation}
\mathbf{E}_{\mathbf{F}}\left[ \left\Vert \mathbf{Z}(l)-\mathbf{Z}%
(l-1)A_{l-1}\right\Vert ^{\kappa }\left\vert \mathcal{B}_{\zeta }\right. %
\right] \leq \left( \mathbf{E}_{\mathbf{F}}\left[ \left\Vert \mathbf{Z}(l)-%
\mathbf{Z}(l-1)A_{l-1}\right\Vert ^{2}\left\vert \,\mathcal{B}_{\zeta
}\right. \right] \right) ^{\kappa /2}.  \label{Nstar}
\end{equation}%
Further we have
\begin{eqnarray}
&&\mathbf{E}_{\mathbf{F}}\left[ \left\Vert \mathbf{Z}(l)-\mathbf{Z}%
(l-1)A_{l-1}\right\Vert ^{2}\left\vert \,\mathcal{B}_{\zeta }\right. \right]
\notag \\
&&\quad =\mathbf{E}_{\mathbf{F}}\left[ \left(
\sum_{i=1}^{m}\sum_{k=1}^{Z_{i}(l-1)}\sum_{j=1}^{m}\left[ \mathbf{\xi }%
_{ij}(l-1;k)-a_{ij}(l-1)\right] \right) ^{2}\left\vert \,\mathcal{B}_{\zeta
}\right. \right]   \notag \\
&&\quad \,=\sum_{i=1}^{m}\mathbf{E}_{\mathbf{F}}\beta _{i}^{2}(l-1)\mathbf{E}%
_{\mathbf{F}}\left[ Z_{i}(l-1)\left\vert \,\mathcal{B}_{\zeta }\right. %
\right]   \notag \\
&=&\sum_{i=1}^{m}\mathbf{E}_{\mathbf{F}}\beta _{i}^{2}(l-1)\left( \mathbf{Z}%
(\zeta )\Pi _{\zeta ,l-1}\right) _{i}\leq \left\Vert \mathbf{Z}(\zeta )\Pi
_{\zeta ,l-1}\right\Vert \sum_{i=1}^{m}\mathbf{E}_{\mathbf{F}}\beta
_{i}^{2}(l-1).  \notag \\
&&\quad \,\,  \label{BB11}
\end{eqnarray}%
Thus, for $\kappa \leq 1$
\begin{equation}
\left( \mathbf{E}_{\mathbf{F}}\left[ \left\Vert \mathbf{Z}(l)-\mathbf{Z}%
(l-1)A_{l-1}\right\Vert ^{2}\left\vert \mathcal{B}_{\zeta }\right. \right]
\right) ^{\kappa /2}\leq \left\Vert \mathbf{Z}(\zeta )\Pi _{\zeta
,l-1}\right\Vert ^{\kappa /2}\left( \sum_{i=1}^{m}\mathbf{E}_{\mathbf{F}%
}\beta _{i}^{2}(l-1)\right) ^{\kappa /2}.  \label{NNstar}
\end{equation}%
This, in view of the inequality
\begin{equation*}
s(\kappa /2)=\lim_{n\rightarrow \infty }\left( \mathbf{E}\left\Vert \Pi
_{0,n}\right\Vert ^{\kappa /2}\right) ^{1/n}<1,
\end{equation*}%
the first part of condition (\ref{Ckappasmall}), and relations (\ref{BBB1})-(%
\ref{NNstar}) leads to the estimate%
\begin{eqnarray}
&&\mathbf{P}\left( \left\vert S(\zeta )-\left\langle \mathbf{Z}(\zeta ),\Xi
_{\zeta }\right\rangle \right\vert >\varepsilon y;\zeta <\infty \right)
\notag \\
&\leq &\frac{K}{\varepsilon ^{\kappa }y^{\kappa }}\mathbf{E}\left[
\sum_{l=\zeta +1}^{\infty }\left( l-\zeta \right) ^{2\kappa }\left\Vert
\mathbf{Z}(\zeta )\right\Vert ^{\kappa /2}\left\Vert \Pi _{\zeta
,l-1}\right\Vert ^{\kappa /2}I\left\{ \zeta <\infty \right\} \right]   \notag
\\
&=&\frac{K}{\varepsilon ^{\kappa }y^{\kappa }}\mathbf{E}\left[ \left\Vert
\mathbf{Z}(\zeta )\right\Vert ^{\kappa /2}I\left\{ \zeta <\infty \right\}
\sum_{l=1}^{\infty }l^{2\kappa }\mathbf{E}\left\Vert \Pi _{0,l-1}\right\Vert
^{\kappa /2}\right]   \notag \\
&\leq &\frac{const}{\varepsilon ^{\kappa }y^{\kappa }r^{\kappa /2}}\mathbf{E}%
\left[ \left\Vert \mathbf{Z}(\zeta )\right\Vert ^{\kappa }I\left\{ \zeta
<\infty \right\} \right] \leq \frac{\varepsilon }{y^{\kappa }}\mathbf{E}%
\left[ \left\Vert \mathbf{Z}(\zeta )\right\Vert ^{\kappa }I\left\{ \zeta
<\infty \right\} \right]   \label{BB2}
\end{eqnarray}%
for all $r\geq r_{0}(\varepsilon )$, proving the lemma for $\kappa
\leq 1 $.

For the case $\kappa >1$ we use Lemma \ref{LGut} to conclude that for any $%
l>\zeta $%
\begin{eqnarray*}
&&\mathbf{E}_{\mathbf{F}}\left[ \left\Vert \mathbf{Z}(l)-\mathbf{Z}%
(l-1)A_{l-1}\right\Vert ^{\kappa }\left\vert \mathcal{B}_{\zeta }\right. %
\right]  \\
&&\quad \leq R_{\kappa }m^{\kappa }\sum_{i=1}^{m}M_{\kappa }(n;i)\mathbf{E}_{%
\mathbf{F}}\left[ \left\Vert Z_{i}(l-1)\right\Vert ^{\kappa /2\vee
1}\left\vert \mathcal{B}_{\zeta }\right. \right]  \\
&&\quad \leq R_{\kappa }m^{\kappa }\mathbf{E}_{\mathbf{F}}\left[ \left\Vert
\mathbf{Z}(l-1)\right\Vert ^{\kappa /2\vee 1}\left\vert \,\mathcal{B}_{\zeta
}\right. \right] \sum_{i=1}^{m}M_{\kappa }(n;i).
\end{eqnarray*}%
By Lemma \ref{LGeom} there exist constants $\rho _{\kappa /2\vee 1}\in
\left( 0,1\right) $ and $K_{\kappa /2\vee 1}<\infty $ such that for all $%
l>\zeta $%
\begin{equation*}
\mathbf{E}_{\mathbf{F}}\left[ \left\Vert
\mathbf{Z}(l-1)\right\Vert ^{\kappa /2\vee 1}\left\vert
\,\mathcal{B}_{\zeta }\right. \right] \leq K_{\kappa /2\vee 1}\rho
_{\kappa /2\vee 1}^{l-\zeta -1}\left\Vert \mathbf{Z}(\zeta
)\right\Vert ^{\kappa /2\vee 1}.
\end{equation*}%
This yields the estimates
\begin{eqnarray*}
&&\mathbf{P}\left( \left\vert S(\zeta )-\left\langle \mathbf{Z}(\zeta ),\Xi
_{\zeta }\right\rangle \right\vert >\varepsilon y;\zeta <\infty \right)  \\
&\leq &\frac{m^{\kappa }K}{\varepsilon ^{\kappa }y^{\kappa }}\mathbf{E}\left[
\sum_{l=\zeta +1}^{\infty }\left( l-\zeta \right) ^{2\kappa }\mathbf{E}_{%
\mathbf{F}}\left[ \left\Vert \mathbf{Z}(l-1)\right\Vert ^{\kappa /2\vee
1}\sum_{i=1}^{m}M_{\kappa }(n;i)\left\vert \mathcal{B}_{\zeta }\right. %
\right] I\left\{ \zeta <\infty \right\} \right]  \\
&\leq &\frac{m^{\kappa +1}KM_{\kappa }}{\varepsilon ^{\kappa }y^{\kappa }}%
\mathbf{E}\left[ \sum_{l=\zeta +1}^{\infty }\left( l-\zeta \right) ^{2\kappa
}\mathbf{E}_{\mathbf{F}}\left[ \left\Vert \mathbf{Z}(l-1)\right\Vert
^{\kappa /2\vee 1}\left\vert \,\mathcal{B}_{\zeta }\right. \right] I\left\{
\zeta <\infty \right\} \right]  \\
&\leq &\frac{m^{\kappa +1}KM_{\kappa }}{\varepsilon ^{\kappa }y^{\kappa }}%
\mathbf{E}\left[ \sum_{l=\zeta +1}^{\infty }\left( l-\zeta \right) ^{2\kappa
}\left\Vert \mathbf{Z}(\zeta )\right\Vert ^{\kappa /2\vee 1}K_{\kappa /2\vee
1}\rho _{\kappa /2\vee 1}^{l-\zeta -1}I\left\{ \zeta <\infty \right\} \right]
\\
&=&\frac{m^{\kappa +1}KM_{\kappa }K_{\kappa /2\vee 1}}{\varepsilon ^{\kappa
}y^{\kappa }}\mathbf{E}\left[ \left\Vert \mathbf{Z}(\zeta )\right\Vert
^{\kappa /2\vee 1}I\left\{ \zeta <\infty \right\} \sum_{l=1}^{\infty
}l^{2\kappa }\rho _{\kappa /2\vee 1}^{l-1}\right]  \\
&\leq &\frac{const}{\varepsilon ^{\kappa }y^{\kappa }r^{\kappa -\kappa
/2\vee 1}}\mathbf{E}\left[ \left\Vert \mathbf{Z}(\zeta )\right\Vert ^{\kappa
}I\left\{ \zeta <\infty \right\} \right] \leq \frac{\varepsilon }{y^{\kappa }%
}\mathbf{E}\left[ \left\Vert \mathbf{Z}(\zeta )\right\Vert ^{\kappa
}I\left\{ \zeta <\infty \right\} \right]
\end{eqnarray*}%
(the last is valid by selecting $r$ sufficiently large) which justifies the
statement of the lemma for $\kappa >1.$

The lemma is proved.

\section{The accumulated amount of the final product}

In this section we deduce some estimates related with the total
size of the final product accumulated in a subcritical MBPRCRE
during its evolution.

Let

\begin{equation*}
\Delta (n):=\sum_{i=1}^{m}\sum_{k=1}^{Z_{i}(n)}\varphi _{i}(n;k)
\end{equation*}%
be the total size of the final product produced by the individuals of the $n$%
-th generation of a MBPRCRE and
\begin{equation*}
\bar{\Phi}(N):=\sum_{n=N}^{\infty }\mathbf{\Delta (}n\mathbf{),\qquad }\Phi =%
\bar{\Phi}_{0}.
\end{equation*}%

\begin{lemma}
\label{Lke222}Under the conditions of Theorem \ref{TailMain} for any $%
\varepsilon >0$ there exists $r=r(\varepsilon )$ such that for all $y\geq
y_{0}$%
\begin{equation}
\mathbf{P}\left( \left\vert \bar{\Phi}(\zeta )-S(\zeta )\right\vert
>\varepsilon y;\zeta <\infty \right) \leq \frac{\varepsilon }{y^{\kappa }}%
\mathbf{E}\left[ \left\Vert \mathbf{Z}(\zeta )\right\Vert ^{\kappa }I\left\{
\zeta <\infty \right\} \right] .  \label{Dif2}
\end{equation}
\end{lemma}

\textbf{Proof}. We have%
\begin{equation*}
\Delta \mathbf{(}n\mathbf{)-}\left\langle \mathbf{Z}(n),\mathbf{C}%
_{n}\right\rangle \mathbf{=}\sum_{i=1}^{m}\sum_{k=1}^{Z_{i}(n)}\left(
\varphi _{i}(n;k)-\mathbf{E}_{\mathbf{F}}\varphi _{i}(n)\right) .
\end{equation*}%
Now we consider separately the cases $\kappa \leq 1$ and $\kappa >1$.

For $\kappa \geq 1$ we use Lemmas \ref{LGut} and \ref{LGeom} to get for $%
n>\zeta :$
\begin{eqnarray*}
&&\mathbf{P}_{\mathbf{F}}\left( \left\vert \Delta \mathbf{(}n\mathbf{)-}%
\left\langle \mathbf{Z}(n),\mathbf{C}_{n}\right\rangle \right\vert >\frac{%
\varepsilon y}{2(n-\zeta )^{2}}\left\vert \mathcal{B}_{\zeta }\right. \right)
\\
&\leq &\frac{4^{\kappa }(n-\zeta )^{2\kappa }}{\left( \varepsilon y\right)
^{\kappa }}\mathbf{E}_{\mathbf{F}}\left[ \left\vert
\sum_{i=1}^{m}\sum_{k=1}^{Z_{i}(n)}\left( \varphi _{i}(n;k)-\mathbf{E}_{%
\mathbf{F}}\varphi _{i}(n)\right) \right\vert ^{\kappa }\left\vert \mathcal{B%
}_{\zeta }\right. \right] \\
&\leq &(4m)^{\kappa }R_{\kappa }\frac{(n-\zeta )^{2\kappa }}{\left(
\varepsilon y\right) ^{\kappa }}\sum_{i=1}^{m}\mathbf{E}_{\mathbf{F}%
}\left\vert \varphi _{i}(n)-\mathbf{E}_{\mathbf{F}}\varphi
_{i}(n)\right\vert ^{\kappa }\mathbf{E}_{\mathbf{F}}\left[ \left\Vert
\mathbf{Z}(n)\right\Vert ^{\kappa /2\vee 1}\left\vert \,\mathcal{B}_{\zeta
}\right. \right] \\
&\leq &(4m)^{\kappa }R_{\kappa }K_{\kappa /2\vee 1}\frac{(n-\zeta )^{2\kappa
}}{\left( \varepsilon y\right) ^{\kappa }}\sum_{i=1}^{m}\mathbf{E}_{\mathbf{F%
}}\left\vert \varphi _{i}(n)-\mathbf{E}_{\mathbf{F}}\varphi
_{i}(n)\right\vert ^{\kappa }\left\Vert \mathbf{Z}(\zeta )\right\Vert
^{\kappa /2\vee 1}\rho _{\kappa /2\vee 1}^{n-\zeta }.
\end{eqnarray*}%
Hence, in view of condition (\ref{Ckappabig})%
\begin{eqnarray*}
&&\mathbf{P}\left( \left\vert \bar{\Phi}(\zeta )-S(\zeta )\right\vert
>\varepsilon y;\zeta <\infty \right) \\
&\leq &\mathbf{E}\left[ \sum_{n=\zeta }^{\infty }\mathbf{P}_{\mathbf{F}%
}\left( \left\vert \Delta \mathbf{(}n\mathbf{)-}\left\langle \mathbf{Z}(n),%
\mathbf{C}_{n}\right\rangle \right\vert >\frac{\varepsilon y}{2(n-\zeta )^{2}%
}\left\vert \,\mathcal{B}_{\zeta }\right. \right) I\left( \zeta <\infty
\right) \right] \\
&\leq &\frac{(4m)^{\kappa }R_{\kappa }K_{\kappa /2\vee 1}K}{\left(
\varepsilon y\right) ^{\kappa }}\mathbf{E}\left[ \left\Vert \mathbf{Z}(\zeta
)\right\Vert ^{\kappa /2\vee 1}\sum_{n=\zeta }^{\infty }(n-\zeta )^{2\kappa
}\rho _{\kappa /2\vee 1}^{n-\zeta }I\left( \zeta <\infty \right) \right] \\
&=&\frac{const}{\left( \varepsilon y\right) ^{\kappa }}\mathbf{E}\left[
\left\Vert \mathbf{Z}(\zeta )\right\Vert ^{\kappa /2\vee 1}I\left( \zeta
<\infty \right) \right] \sum_{n=1}^{\infty }n^{2\kappa }\rho _{\kappa /2\vee
1}^{n} \\
&\leq &\frac{const}{\left( \varepsilon y\right) ^{\kappa }r^{\kappa -(\kappa
/2\vee 1)}}\mathbf{E}\left[ \left\Vert \mathbf{Z}(\zeta )\right\Vert
^{\kappa }I\left( \zeta <\infty \right) \right] \leq \frac{\varepsilon }{%
y^{\kappa }}\mathbf{E}\left[ \left\Vert \mathbf{Z}(\zeta )\right\Vert
^{\kappa }I\left( \zeta <\infty \right) \right]
\end{eqnarray*}%
for all $r\geq r(\varepsilon ).$

To analyze the case $\kappa \leq 1$ we apply for $n>\zeta $ the inequality%
\begin{equation*}
\mathbf{E}_{\mathbf{F}}\left[ \left\vert \Delta \mathbf{(}n\mathbf{)-}%
\left\langle \mathbf{Z}(n),\mathbf{C}_{n}\right\rangle \right\vert ^{\kappa
}\left\vert \mathcal{B}_{\zeta }\right. \right] \leq \left( \mathbf{E}_{%
\mathbf{F}}\left[ \left\vert \Delta \mathbf{(}n\mathbf{)-}\left\langle
\mathbf{Z}(n),\mathbf{C}_{n}\right\rangle \right\vert ^{2}\left\vert \,%
\mathcal{B}_{\zeta }\right. \right] \right) ^{\kappa /2}.
\end{equation*}%
Further, we have
\begin{eqnarray*}
&&\mathbf{E}_{\mathbf{F}}\left[ \left\vert \Delta \mathbf{(}n\mathbf{)-}%
\left\langle \mathbf{Z}(n),\mathbf{C}_{n}\right\rangle \right\vert
^{2}\left\vert \mathcal{B}_{\zeta }\right. \right] \leq m^{2}\sum_{i=1}^{m}%
\mathbf{E}_{\mathbf{F}}\beta _{i}^{2}(n)\mathbf{E}_{\mathbf{F}}\left[
Z_{i}(n)\left\vert \,\mathcal{B}_{\zeta }\right. \right] \\
&=&m^{2}\sum_{i=1}^{m}\mathbf{E}_{\mathbf{F}}\beta _{i}^{2}(n)\left( \mathbf{%
Z}(\zeta )\Pi _{\zeta ,n}\right) _{i}\leq m^{2}\left\Vert \mathbf{Z}(\zeta
)\Pi _{\zeta ,n}\right\Vert \sum_{i=1}^{m}\mathbf{E}_{\mathbf{F}}\beta
_{i}^{2}(n).
\end{eqnarray*}%
Thus, for $\kappa \leq 1$
\begin{equation*}
\left( \mathbf{E}_{\mathbf{F}}\left[ \left\vert \Delta \mathbf{(}n\mathbf{)-}%
\left\langle \mathbf{Z}(n),\mathbf{C}_{n}\right\rangle \right\vert
^{2}\left\vert \,\mathcal{B}_{\zeta }\right. \right] \right) ^{\kappa
/2}\leq m^{\kappa }\left\Vert \mathbf{Z}(\zeta )\Pi _{\zeta ,n+1}\right\Vert
^{\kappa /2}\left( \sum_{i=1}^{m}\mathbf{E}_{\mathbf{F}}\beta
_{i}^{2}(n)\right) ^{\kappa /2}.
\end{equation*}%
This combined with the assumption (\ref{Ckappasmall}) shows that%
\begin{eqnarray*}
&&\mathbf{P}\left( \left\vert \bar{\Phi}(\zeta )-S(\zeta )\right\vert
>\varepsilon y;\zeta <\infty \right) \\
&\leq &\mathbf{E}\left[ \sum_{n=\zeta }^{\infty }\mathbf{P}_{\mathbf{F}%
}\left( \left\vert \Delta \mathbf{(}n\mathbf{)-}\left\langle \mathbf{Z}(n),%
\mathbf{C}_{n}\right\rangle \right\vert >\frac{\varepsilon y}{2(n-\zeta )^{2}%
}\left\vert \,\zeta ,\mathbf{Z}(0),\ldots ,\mathbf{Z}(\zeta )\right. \right)
I\left( \zeta <\infty \right) \right] \\
&\leq &\frac{m^{\kappa }K}{\varepsilon ^{\kappa }y^{\kappa }}\mathbf{E}\left[
\sum_{l=\zeta +1}^{\infty }\left( l-\zeta \right) ^{2\kappa }\left\Vert
\mathbf{Z}(\zeta )\right\Vert ^{\kappa /2}\left\Vert \Pi _{\zeta
,l}\right\Vert ^{\kappa /2}I\left( \zeta <\infty \right) \right] \\
&=&\frac{m^{\kappa }K}{\varepsilon ^{\kappa }y^{\kappa }}\mathbf{E}\left[
\left\Vert \mathbf{Z}(\zeta )\right\Vert ^{\kappa /2}\sum_{l=1}^{\infty
}l^{2\kappa }\mathbf{E}\left\Vert \Pi _{0,l}\right\Vert ^{\kappa /2}I\left(
\zeta <\infty \right) \right] \\
&\leq &\frac{const}{\varepsilon ^{\kappa }y^{\kappa }r^{\kappa /2}}\mathbf{E}%
\left[ \left\Vert \mathbf{Z}(\zeta )\right\Vert ^{\kappa }I\left\{ \zeta
<\infty \right\} \right] \leq \frac{\varepsilon }{y^{\kappa }}\mathbf{E}%
\left[ \left\Vert \mathbf{Z}(\zeta )\right\Vert ^{\kappa }I\left\{ \zeta
<\infty \right\} \right]
\end{eqnarray*}%
for all $r\geq r(\varepsilon )$.

The lemma is proved.

Up to now we have assumed that the initial number of particles and
the initial
volume of the final product are nonrandom. Lemmas \ref{Lremaind} and \ref%
{LatInfty} are free of this restriction.

\begin{lemma}
\label{Lremaind}Let a MBPRCRE be subcritical, $\mathbf{E}\left\Vert \mathbf{Z%
}(0)\right\Vert <\infty $ and there exist $\delta >0$ such that
\begin{equation*}
\mathbf{E}\Phi ^{\kappa +\delta }(0)<\infty ,\quad \max_{1\leq i\leq m}%
\mathbf{E}\varphi _{i}^{\kappa +\delta }<\infty .
\end{equation*}%
Then for any $\varepsilon >0$ there exists $r=r(\varepsilon )$ such that for
all $y\geq y_{0}$
\begin{equation}
\mathbf{P}\left( \Phi (\zeta -1)>\varepsilon y;\zeta <\infty \right) <\frac{1%
}{y^{\kappa +\delta /2}}.  \label{Dif3}
\end{equation}
\end{lemma}

\textbf{Proof}. Let, as earlier, $\tau $ be the extinction moment of the
underlying MBPRE and let $c$ be a constant such that for all $y\geq y_{0}(c)$
\begin{equation*}
\sum_{t>c\ln y}\mathbf{P}\left( \tau >t\right) \leq \frac{1}{y^{\kappa
+\delta }}.
\end{equation*}%
Such a constant, clearly, exists in view of (\ref{SurUp}). Recalling (\ref%
{DDphi}) put%
\begin{equation*}
\Phi _{r}(n)=\Phi (0)+\sum_{l=0}^{n}\sum_{i=1}^{m}\sum_{k=1}^{r}\varphi
_{i}(l;k).
\end{equation*}%
Now we have for $y\geq y_{0}(c)$
\begin{eqnarray*}
&&\mathbf{P}\left( \Phi (\zeta -1)>\varepsilon y;\zeta <\infty \right)
=\sum_{n=1}^{\infty }\mathbf{P}\left( \Phi (n-1)>\varepsilon y;\zeta
=n\right) \\
&\leq &\sum_{n=1}^{\infty }\mathbf{P}\left( \Phi _{r}(n-1)>\varepsilon
y;\zeta =n\right) \leq \sum_{n=1}^{\infty }\mathbf{P}\left( \Phi
_{r}(n-1)>\varepsilon y;\tau >n\right) \\
&\leq &\mathbf{P}\left( \Phi (0)>\frac{\varepsilon y}{2}\right) +\sum_{1\leq
n\leq c\ln y}\sum_{l=0}^{n}\sum_{i=1}^{m}\sum_{k=1}^{r}\mathbf{P}\left(
\varphi _{i}(l;k)>\frac{\varepsilon y}{2mrc^{2}\ln ^{2}y}\right) \\
&&+\sum_{n>c\ln y}\mathbf{P}\left( \tau >n\right) \leq \frac{2}{\left(
\varepsilon y\right) ^{\kappa +\delta }}\mathbf{E}\Phi ^{\kappa +\delta }(0)
\\
&&+\frac{\left( 2mrc^{2}\ln ^{2}y\right) ^{\kappa +\delta }}{\left(
\varepsilon y\right) ^{\kappa +\delta }}\sum_{1\leq n\leq c\ln
y}\sum_{l=0}^{n}\sum_{i=1}^{m}\sum_{k=1}^{r}\mathbf{E}\varphi _{i}^{\kappa
+\delta } +\frac{1}{y^{\kappa +\delta }} \\
&&\quad \leq K\frac{\ln ^{2(\kappa +\delta +1)}y}{\left( \varepsilon
y\right) ^{\kappa +\delta }}+\frac{1}{y^{\kappa +\delta }}\leq \frac{1}{%
y^{\kappa +\delta /2}}.
\end{eqnarray*}

\begin{lemma}
\label{LatInfty}Under conditions of Lemma \ref{Lremaind} for any $r$ there
exists $y_{0}=y_{0}(r)$ such that for all $y>y_{0}$
\begin{equation}
\mathbf{P}\left( \Phi >y;\zeta =\infty \right) \leq \frac{1}{y^{\kappa
+\delta /2}}.  \label{Dif4}
\end{equation}
\end{lemma}

\textbf{Proof}. Letting $c$ be the same as in the previous lemma we have for
$y\geq y_{0}(c)$%
\begin{eqnarray*}
&&\mathbf{P}\left( \Phi >y;\zeta =\infty \right) \leq \mathbf{P}\left( \Phi
_{r}(\tau )>y\right) \leq \mathbf{P}\left( \Phi _{r}(\left[ c\ln y\right]
)>y\right) +\mathbf{P}\left( \tau >c\ln y\right) \\
&\leq &\mathbf{P}\left( \Phi (0)>\frac{y}{2}\right) +\sum_{l=0}^{\left[ c\ln
y\right] }\sum_{i=1}^{m}\sum_{k=1}^{r}\mathbf{P}\left( \varphi _{i}(l;k)>%
\frac{y}{2mr\left[ c\ln y\right] }\right) \\
&&+\mathbf{P}\left( \tau >c\ln y\right) \leq K\frac{\ln ^{\kappa +\delta +1}y%
}{y^{\kappa +\delta }}+\frac{1}{y^{\kappa +\delta }}
\end{eqnarray*}%
as desired.

\section{Proof of Theorem \protect\ref{TailMain}\label{MMain}}

Now we are ready to prove the main result of the paper, Theorem \ref%
{TailMain}. First observe that by the equivalence of the norms $\left\Vert
\cdot \right\Vert$ and $\left\Vert \cdot\right\Vert_2$ and estimates (\ref%
{Diff1}),(\ref{Dif2}), (\ref{Dif3}) and (\ref{Dif4}), for any
$\varepsilon \in (0,1/3)$ one can find $r=r(\varepsilon )$ such
that for all $y\geq y_{0}(r,\varepsilon )$
\begin{eqnarray}
&&\mathbf{P}\left( \Phi >y\right) \leq \mathbf{P}\left( \Phi >y;\zeta
=\infty \right)  \notag \\
&& +\mathbf{P}\left( \left\langle \mathbf{Z}(\zeta ),\Xi _{\zeta
}\right\rangle >y(1-3\varepsilon );\zeta <\infty \right)+\mathbf{P}\left(
\left\vert \bar{\Phi}(\zeta )-S(\zeta )\right\vert >\varepsilon y;\zeta
<\infty \right)  \notag \\
&&+\mathbf{P}\left( \left\vert S(\zeta )-\left\langle \mathbf{Z}(\zeta ),\Xi
_{\zeta }\right\rangle \right\vert >\varepsilon y;\zeta <\infty \right) +%
\mathbf{P}\left( \Phi (\zeta -1)>\varepsilon y;\zeta <\infty \right)  \notag
\\
&\leq &\mathbf{P}\left( \left\langle \mathbf{Z}(\zeta ),\Xi _{\zeta
}\right\rangle >y(1-3\varepsilon );\zeta <\infty \right) +\frac{2\varepsilon
}{y^{\kappa }}\mathbf{E}\left[ \left\Vert \mathbf{Z}(\zeta )\right\Vert_2
^{\kappa }I\left\{ \zeta <\infty \right\} \right] +\frac{2}{y^{\kappa
+\delta /2}}.  \notag \\
&&  \label{FiDecomp}
\end{eqnarray}%
Let $K(l):=\inf_{\mathbf{u\in U}_{+}}l(\mathbf{u})>0$ where $l(\mathbf{u})$
is the function involved in Condition~$T$. By this condition and the
independency of $\mathbf{Z}(\zeta )$ and $\Xi _{\zeta }$ we conclude
\begin{eqnarray*}
&&\lim \sup_{y\rightarrow \infty }y^{\kappa }\mathbf{P}\left( \left\langle
\mathbf{Z}(\zeta ),\Xi _{\zeta }\right\rangle >y(1-3\varepsilon );\zeta
<\infty \right) \\
&\leq &\lim \sup_{y\rightarrow \infty }y^{\kappa }\int_{\left\Vert \mathbf{u}%
\right\Vert_2 =r}^{\infty }\mathbf{P}\left( \mathbf{Z}(\zeta )\in d\mathbf{u}%
;\zeta <\infty \right) \mathbf{P}\left( \left\langle \frac{\mathbf{u}}{%
\left\Vert \mathbf{u}\right\Vert_2 },\Xi \right\rangle \geq \frac{%
(1-3\varepsilon )y}{\left\Vert \mathbf{u}\right\Vert_2 }\right) \\
&=&K_{0}(1-3\varepsilon )^{-\kappa }\int_{\left\Vert \mathbf{u}\right\Vert_2
=r}^{\infty }\mathbf{P}\left( \mathbf{Z}(\zeta )\in d\mathbf{u};\zeta
<\infty \right) \left\Vert \mathbf{u}\right\Vert_2^{\kappa }l\left( \frac{%
\mathbf{u}}{\left\Vert \mathbf{u}\right\Vert_2}\right) \\
&=&K_{0}(1-3\varepsilon )^{-\kappa }\mathbf{E}\left[ \left\Vert \mathbf{Z}%
(\zeta )\right\Vert_2 ^{\kappa }l\left( \frac{\mathbf{Z}(\zeta )}{\left\Vert
\mathbf{Z}(\zeta )\right\Vert_2}\right) I\left\{\zeta <\infty \right\} %
\right] <\infty .
\end{eqnarray*}%
Since $y^{\kappa }\mathbf{P}\left( \Phi >y\right) $ does not
depend on $r$ and $\varepsilon, $  the previous estimate and
(\ref{FiDecomp}) yield
\begin{equation}
\lim \sup_{y\rightarrow \infty }y^{\kappa }\mathbf{P}\left( \Phi >y\right)
<\infty  \label{PPinf}
\end{equation}%
and, moreover,%
\begin{equation}
\lim \sup_{y\rightarrow \infty }y^{\kappa }\mathbf{P}\left( \Phi >y\right)
\leq K_{0}\lim_{r\rightarrow \infty }\mathbf{E}\left[ \left\Vert \mathbf{Z}%
(\zeta )\right\Vert_2^{\kappa }l\left( \frac{\mathbf{Z}(\zeta )}{\left\Vert
\mathbf{Z}(\zeta )\right\Vert_2 }\right) I\left\{ \zeta <\infty \right\} %
\right].  \label{PP2}
\end{equation}

To get a similar estimate from below we use for $\varepsilon >0$ the
inequality%
\begin{eqnarray*}
&&\mathbf{P}\left( \Phi >y\right) \geq\mathbf{P}\left( \Phi >y;\zeta <\infty
\right) \geq \mathbf{P}\left( \left\langle \mathbf{Z}(\zeta ),\Xi _{\zeta
}\right\rangle >y(1+3\varepsilon );\zeta <\infty \right) \\
&&\quad- \mathbf{P}\left( \left\vert \bar{\Phi}(\zeta )-S(\zeta )\right\vert
>\varepsilon y;\zeta <\infty \right) -\mathbf{P}\left( \left\vert S(\zeta
)-\left\langle \mathbf{Z}(\zeta ),\Xi _{\zeta }\right\rangle \right\vert
>\varepsilon y;\zeta <\infty \right) \\
&&\quad- \mathbf{P}\left( \Phi \left( \zeta -1\right) >\varepsilon y;\zeta
<\infty \right) .
\end{eqnarray*}%
Now we select $r$ as large to meet estimates
(\ref{Diff1}), (\ref{Dif2}) and (\ref{Dif3}). This gives for sufficiently large $y>r$ the inequality%
\begin{equation*}
\mathbf{P}\left( \Phi >y\right) \geq \mathbf{P}\left( \left\langle \mathbf{Z}%
(\zeta ),\Xi _{\zeta }\right\rangle >y(1+3\varepsilon );\zeta <\infty
\right) -\frac{2\varepsilon }{y^{\kappa }}\mathbf{E}\left[ \left\Vert
\mathbf{Z}(\zeta )\right\Vert_2 ^{\kappa }I\left\{ \zeta <\infty \right\} %
\right] -\frac{1}{y^{\kappa +\delta /2}}.
\end{equation*}%
Letting $y\rightarrow \infty $ we obtain%
\begin{eqnarray*}
&&\lim \inf_{y\rightarrow \infty }y^{\kappa }\mathbf{P}\left( \left\langle
\mathbf{Z}(\zeta ),\Xi _{\zeta }\right\rangle >y(1+3\varepsilon )\right) \\
&=&\lim \inf_{y\rightarrow \infty }y^{\kappa }\int_{\left\Vert \mathbf{u}%
\right\Vert_2 =r}^{\infty }\mathbf{P}\left( \mathbf{Z}(\zeta )\in d\mathbf{u}%
;\zeta <\infty \right) \mathbf{P}\left( \left\langle \frac{\mathbf{u}}{%
\left\Vert \mathbf{u}\right\Vert_2 },\Xi \right\rangle \geq \frac{%
(1+3\varepsilon )y}{\left\Vert \mathbf{u}\right\Vert_2 }\right) \\
&=&K_{0}(1+3\varepsilon )^{-\kappa }\mathbf{E}\left[ \left\Vert \mathbf{Z}%
(\zeta )\right\Vert_2^{\kappa }l\left( \frac{\mathbf{Z}(\zeta )}{\left\Vert
\mathbf{Z}(\zeta )\right\Vert_2 }\right) I\left\{ \zeta <\infty \right\} %
\right] .
\end{eqnarray*}%
We know that for sufficiently small $\varepsilon >0$ and an appropriate $r$%
\begin{eqnarray*}
&&K_{0}(1+3\varepsilon )^{-\kappa }\mathbf{E}\left[ \left\Vert \mathbf{Z}%
(\zeta )\right\Vert_2 ^{\kappa }l\left( \frac{\mathbf{Z}(\zeta )}{\left\Vert
\mathbf{Z}(\zeta )\right\Vert_2 }\right) I\left\{ \zeta <\infty \right\} %
\right] -2\varepsilon \mathbf{E}\left[ \left\Vert \mathbf{Z}(\zeta
)\right\Vert_2 ^{\kappa }I\left\{ \zeta <\infty \right\} \right] \\
&\geq &\left( K_{0}K(l)(1+3\varepsilon )^{-\kappa }-2\varepsilon \right)
\mathbf{E}\left[ \left\Vert \mathbf{Z}(\zeta )\right\Vert_2 ^{\kappa
}I\left\{ \zeta <\infty \right\} \right] >0
\end{eqnarray*}%
which implies
\begin{equation*}
\lim \inf_{y\rightarrow \infty }y^{\kappa }\mathbf{P}\left( \Phi >y\right) >0
\end{equation*}%
leading in turn to%
\begin{equation}
\lim \inf_{y\rightarrow \infty }y^{\kappa }\mathbf{P}\left( \Phi >y\right)
\geq K_{0}\lim_{r\rightarrow \infty }\mathbf{E}\left[ \left\Vert \mathbf{Z}%
(\zeta )\right\Vert_2 ^{\kappa }l\left( \frac{\mathbf{Z}(\zeta )}{\left\Vert
\mathbf{Z}(\zeta )\right\Vert_2 }\right) I\left\{ \zeta <\infty \right\} %
\right] .  \label{FiBelow}
\end{equation}%
This combined with (\ref{PP2}) gives
\begin{equation*}
\lim_{y\rightarrow \infty }y^{\kappa }\mathbf{P}\left( \Phi >y\right)
=K_{0}\lim_{r\rightarrow \infty }\mathbf{E}\left[ \left\Vert \mathbf{Z}%
(\zeta )\right\Vert_2 ^{\kappa }l\left( \frac{\mathbf{Z}(\zeta )}{\left\Vert
\mathbf{Z}(\zeta )\right\Vert_2 }\right) I\left\{ \zeta <\infty \right\} %
\right] \in (0,\infty ).
\end{equation*}

The theorem is proved.

\begin{theorem}
\label{TfinalZero}Let conditions of Theorem \ref{TailMain} be valid for a
subcritical MBPRCRE starting at moment $0$ by a random tuple $\mathbf{Z}%
(0)=(Z_{1}(0),\ldots,Z_{m}(0))$ of particles with $\mathbf{P}\left( \mathbf{Z%
}(0)\neq \mathbf{0}\right) >0$ and having a random initial size $\Phi (0)$
of the final product. If
\begin{equation*}
\mathbf{E}\left[ \left\Vert \mathbf{Z}(0)\right\Vert ^{t\vee 1}+\Phi ^{t}(0)%
\right] <\infty
\end{equation*}%
for any $t\in \left( 0,\kappa \right) $, then
\begin{equation*}
\mathbf{E}\Phi ^{x}<\infty
\end{equation*}%
if and only if $x\in \left( 0,\kappa \right) $.
\end{theorem}

\textbf{Proof}. Let $\Phi _{ik}$ be the total size of the final
product produced by all descendants of the $k$-th particle of type
$i$ of the zero generation. Clearly, the accumulated amount $\Phi
$ of the final product can be written as
\begin{equation*}
\Phi =\Phi (0)+\sum_{i=1}^{m}\sum_{k=1}^{Z_{i}(0)}\Phi _{ik}.
\end{equation*}%
Given $\mathbf{Z}(0)=\mathbf{z}=(z_{1},\ldots ,z_{m})$ we have for $x\leq 1$%
\begin{eqnarray*}
\mathbf{E}\left[ \Phi ^{x}|\mathbf{Z}(0)=\mathbf{z}\right] &\leq &\mathbf{E}%
\Phi ^{x}(0)+\sum_{i=1}^{m}\sum_{k=1}^{z_{i}}\mathbf{E}\Phi _{ik}^{x}=%
\mathbf{E}\Phi ^{x}(0)+\sum_{i=1}^{m}z_{i}\mathbf{E}\left[ \Phi ^{x}|\mathbf{%
Z}(0)=\mathbf{e}_{i}\right] \\
&\leq &\mathbf{E}\Phi ^{x}(0)+\left\Vert \mathbf{z}\right\Vert \max_{1\leq
i\leq m}\mathbf{E}\left[ \Phi ^{x}|\mathbf{Z}(0)=\mathbf{e}_{i}\right] ,
\end{eqnarray*}%
while for $x>1$ there exists a constant $R_{x}^{\ast }$ such that (see, for
instance, Theorem~5.2, page 22 in \cite{Gut})
\begin{eqnarray*}
\mathbf{E}\left[ \Phi ^{x}|\mathbf{Z}(0)=\mathbf{z}\right] &\leq &\left(
m+1\right) ^{x}\left[ \mathbf{E}\Phi ^{x}(0)+\sum_{i=1}^{m}\mathbf{E}\left(
\sum_{k=1}^{z_{i}}\Phi _{ik}\right) ^{x}\right] \\
&\leq &\left( m+1\right) ^{x}\left[ \mathbf{E}\Phi ^{x}(0)+R_{x}^{\ast
}\sum_{i=1}^{m}z_{i}^{x}\mathbf{E}\left[ \Phi ^{x}|\mathbf{Z}(0)=\mathbf{e}%
_{i}\right] \right] \\
&\leq &\left( m+1\right) ^{x}\left[ \mathbf{E}\Phi ^{x}(0)+R_{x}^{\ast
}\left\Vert \mathbf{z}\right\Vert ^{x}\max_{1\leq i\leq m}\mathbf{E}\left[
\Phi ^{x}|\mathbf{Z}(0)=\mathbf{e}_{i}\right] \right] .
\end{eqnarray*}%
By the total probability formula and the estimates above we obtain for $x\in
(0,\kappa )$:
\begin{eqnarray*}
\mathbf{E}\Phi ^{x} &=&\mathbf{P}\left( \mathbf{Z}(0)=\mathbf{0}\right)
\mathbf{E}\Phi ^{x}(0)+\sum_{\mathbf{z}\in \mathbb{Z}_{+}\backslash \left\{
0\right\} }\mathbf{P}\left( \mathbf{Z}(0)=\mathbf{z}\right) \mathbf{E}\left[
\Phi ^{x}|\mathbf{Z}(0)=\mathbf{z}\right] \\
&\leq &K_{2}\mathbf{E}\Phi ^{x}(0)+K_{3}\max_{1\leq i\leq m}\mathbf{E}\left[
\Phi ^{x}|\mathbf{Z}(0)=\mathbf{e}_{i}\right] \mathbf{E}\left[ \left\Vert
\mathbf{Z}(0)\right\Vert ^{x\vee 1}\right] <\infty .
\end{eqnarray*}%
On the other hand, for any $\mathbf{z}\in
\mathbb{N}_{0}^m\backslash \left\{ \mathbf{0}\right\} $ such that
$\mathbf{P}\left( \mathbf{Z}(0)=\mathbf{z}\right) >0$
\begin{equation*}
\mathbf{E}\Phi ^{x}\geq \mathbf{E}\left[ \Phi ^{x}|\mathbf{Z}(0)=\mathbf{z}%
\right] \mathbf{P}\left( \mathbf{Z}(0)=\mathbf{z}\right)
\end{equation*}%
and the desired result for $x\geq \kappa $ follows from Theorem \ref%
{TailMain}.

\section{Polling systems with zero switchover times and MBPRCRE\label%
{SecPoll}}

We start this section by the description of a connection between the
BTPSFPRE with zero switchover times and the MBPRCRE.

Consider a polling system with zero switchover times and assume that the
operation of the initially idle system starts at the moment when a customer
arrives to a station $J\in \left\{ 1,2,\ldots,m\right\} $. We would like to
study the distribution of the busy period of the server which performs $J-1$%
\, switches of zero length and then starts the service of the
customer arrived to station $J$.

The length of this busy period is constituted by the time intervals spend by
the server to make a random number of complete cycles $(1\rightarrow
2\rightarrow \cdots \rightarrow m\rightarrow 1)$ plus the time needed to
perform the last (may be, incomplete) route just before the moment when
there are no customers in the system for the first time.

Assume that in the course of the $n-$th service cycle within the
initial busy period the service discipline at station $i$
satisfies the branching property with vector-valued m.p.g.f.
\begin{equation*}
\mathbf{\phi }_{n}(\mathbf{s;}\lambda )=\left( \phi _{n}^{(1)}(\mathbf{s;}%
\lambda ),\ldots ,\phi _{n}^{(m)}(\mathbf{s;}\lambda )\right)
\end{equation*}%
where
\begin{equation*}
\phi _{n}^{(i)}(\mathbf{s;}\lambda ):=\mathbf{E}\left[ s_{1}^{\theta
_{i1}(n)}s_{2}^{\theta _{i2}(n)}\cdots s_{m}^{\theta _{im}(n)}e^{-\lambda
\phi _{i}(n)}\right] ,\,i=1,2,\ldots ,m.
\end{equation*}%
Suppose that the sequence $\mathbf{\phi }_{0}(\mathbf{s;}\lambda ),\mathbf{%
\phi }_{1}(\mathbf{s;}\lambda ),\ldots $ is selected at random in an iid
manner. Let further, for each customer, say $j$, served at station $i$
during the $n-$th route the final product $\phi _{i}(n,j)$ and the vector $%
\left( \theta _{i1}(n,j),\ldots ,\theta _{im}(n,j)\right) $ of the numbers
of new customers arrived to the system during the service time $\tau
_{i}(n,j)$ of the customer under consideration have the property%
\begin{equation*}
\left( \theta _{i1}(n,j),\ldots ,\theta _{im}(n,j);\phi _{i}(n,j)\right)
\overset{d}{=}\left( \theta _{i1}(n),\ldots ,\theta _{im}(n);\phi
_{i}(n)\right).
\end{equation*}%
Here the final product may be not only $\tau _{i}(n,j)$ but any nonnegative
random variable being either dependent on the the tuple $\left( \theta
_{i1}(n,j),\ldots ,\theta _{im}(n,j);\tau _{i}(n,j)\right) $ or independent
on the performance of the system at all. Set $h_{n}^{(i)}(\mathbf{s}):=\phi
_{n}^{(i)}(\mathbf{s};0),i=1,2,\ldots ,m$ and for $n=0,1,2,\ldots $
introduce m.p.g.f.'s
\begin{equation*}
F_{n}^{(i)}(\mathbf{s;\lambda })=\mathbf{E}\left[ s_{1}^{\xi
_{i1}(n)}s_{2}^{\xi _{i2}(n)}\cdots s_{m}^{\xi _{im}(n)}e^{-\lambda \varphi
_{i}(n)}\right]
\end{equation*}%
and p.g.f.'s
\begin{equation*}
f_{n}^{(i)}(\mathbf{s})=\mathbf{E}\left[ s_{1}^{\xi _{i1}(n)}s_{2}^{\xi
_{i2}(n)}\cdots s_{m}^{\xi _{im}(n)}\right]
\end{equation*}%
by the equalities $F_{n}^{(m)}(\mathbf{s};\lambda )=\phi _{n}^{(m)}(\mathbf{%
s;}\lambda )$,
\begin{equation}
F_{n}^{(i)}(\mathbf{s};\lambda )=\phi _{n}^{(i)}\left( s_{1},\ldots
,s_{i},F_{n}^{(i+1)}\left( \mathbf{s};\lambda \right) ,\ldots
,F_{n}^{(m)}\left( \mathbf{s};\lambda \right) ;\lambda \right) ,\,i<m,
\label{GG3}
\end{equation}%
and $f_{n}^{(m)}(\mathbf{s})=h_{n}^{(m)}(\mathbf{s}),$
\begin{equation}
f_{n}^{(i)}(\mathbf{s})=h_{n}^{(i)}\left( s_{1},\ldots
,s_{i},f_{n}^{(i+1)}\left( \mathbf{s}\right) ,\ldots ,f_{n}^{(m)}\left(
\mathbf{s}\right) \right) ,\,i<m.  \label{GGGG}
\end{equation}

We would like to describe conditions on the branching type polling system
under which power moments of the amount of the final product accumulated in
the system during its busy period are finite or infinite. Our results are
based on the following a bit long but important statement revealing
connections between the behavior of certain characteristics of the busy
periods of BTPSFPRE and related characteristics of MBPRCRE whose population
and the size of the final product at moment $0$ are random.

\begin{theorem}
\label{ResRep}The joint distribution of the number of customers at different
stations at the end of the $n-$th service cycle and the amount of the final
product accumulated in the system to the end of the $n-$th service cycle in
a BTPSFPRE starting by a single customer at station $J$ and stopped at the
end of the first busy period coincides with the joint distribution of the
number of particles in the $n-$th generation and the total amount of the
final product produced for the $n-1$ generations in a MBPRCRE whose
generation $0$ is specified by a random number of particles and a random
amount of the final product with m.p.g.f. $F_{0}^{(J)}(\mathbf{s;}\lambda )$
and where the joint distribution of the number of direct descendants and the
amount of the final product produced by particles of different types of the $%
k$-th generation is given by the vector-valued m.p.g.f.'s
\begin{equation*}
\mathbf{F}_{k}(\mathbf{s;}\lambda )=\left( F_{k}^{(1)}(\mathbf{s;}\lambda
),\ldots ,F_{k}^{(m)}(\mathbf{s;}\lambda )\right) ,\,k=1,2,\ldots ,n.
\end{equation*}
\end{theorem}

To prove this theorem one should repeat almost literally the proof of
Theorem 4 in \cite{Res93} and we omit the respective arguments.

Note that if, instead of a single individual at station $J$ we
would initially have in the system a batch of customers $\left(
k_{1},\ldots ,k_{m}\right) $ with $k_{i}$ customers at station~$i$
then the distribution of the initial number of particles and the
initial size of the final product in the corresponding MBPRCRE
should be specified by the m.p.g.f.
\begin{equation*}
F_{0}(\mathbf{s;}\lambda )=\prod_{J=1}^{m}\left( F_{0}^{(J)}(\mathbf{s;}%
\lambda )\right) ^{k_{J}}.
\end{equation*}

We call the MBPRCRE described by Theorem \ref{ResRep} the \textit{associated}
MBPRCRE for the BTPSFPRE.

Now we may reformulate the results of Section \ref{SecLimT} in terms of our
polling system.

Let $A_{n}:=\left( a_{ij}(n)\right) _{i,j=1}^{m}$ be the matrix with elements%
\begin{equation}
a_{ij}(n):=\frac{\partial }{\partial s_{j}}f_{n}^{(i)}(\mathbf{s})\left\vert
_{\mathbf{s}=\mathbf{1}}\right. =\mathbf{E}_{\mathbf{f}}\xi _{ij}(n)
\label{DefElem}
\end{equation}%
and let $H_{n}:=\left( h_{ij}(n)\right) _{i,j=1}^{m}$ be the matrix%
\begin{equation*}
h_{ij}(n):=\frac{\partial }{\partial s_{j}}h_{n}^{(i)}(\mathbf{s})\left\vert
_{\mathbf{s}=\mathbf{1}}\right. =\mathbf{E}_{\mathbf{h}}\theta _{ij}(n).
\end{equation*}%
Then in view of (\ref{GGGG}) $a_{mj}(n)=h_{mj}(n),j=1,2,\ldots ,m$ and for $%
i<m$
\begin{equation}
a_{ij}(n)=h_{ij}(n)I\left\{ j\leq i\right\}
+\sum_{k=i+1}^{m}h_{ik}(n)a_{kj}(n).  \label{lMatem}
\end{equation}%
For $i=1,\ldots ,m$ introduce auxiliary matrices%
\begin{equation*}
H_{n}^{(i)}:=\left(
\begin{array}{ccccccccc}
1 & 0 & \cdots & 0 & 0 & 0 & \cdots & \cdots & 0 \\
0 & 1 & \overset{\cdot }{}\cdot \underset{\cdot }{} & \overset{\cdot }{%
\underset{\cdot }{\cdot }} & 0 & 0 & \cdots & \cdots & 0 \\
\overset{\cdot }{\underset{\cdot }{\cdot }} & \overset{\cdot }{}\cdot
\underset{\cdot }{} & \overset{\cdot }{}\cdot \underset{\cdot }{} & 0 & 0 & 0
& \cdots & \cdots & 0 \\
0 & \cdots & 0 & 1 & 0 & 0 & \cdots & \cdots & 0 \\
h_{i1}(n) & h_{i2}(n) & \cdots & h_{i(i-1)}(n) & h_{ii}(n) & h_{i(i+1)}(n) &
\overset{\cdot }{\underset{\cdot }{\cdot }} & \overset{\cdot }{\underset{%
\cdot }{\cdot }} & h_{im}(n) \\
0 & \cdots & \cdots & 0 & 0 & 1 & 0 & \cdots & 0 \\
0 & \cdots & \cdots & 0 & 0 & 0 & 1 & \overset{\cdot }{}\cdot \underset{%
\cdot }{} & 0 \\
0 & \cdots & \cdots & 0 & 0 & 0 & \overset{\cdot }{}\cdot \underset{\cdot }{}
& \overset{\cdot }{}\cdot \underset{\cdot }{} & 0 \\
0 & \cdots & \cdots & 0 & 0 & 0 & \cdots & 0 & 1%
\end{array}%
\right) ,
\end{equation*}%
where for each $i$ the elements of the matrix $H_{n}$ are located only in
row $i$ of the matrix $H_{n}^{(i)}$. It is not difficult to check by (\ref%
{lMatem}) that%
\begin{equation}
A_{n}=H_{n}^{(1)}H_{n}^{(2)}\cdots H_{n}^{(m)}.  \label{Matexplic}
\end{equation}

Further, let $\mathbf{C}_{n}:=\left( C_{1}(n),\ldots,C_{m}(n)\right)
^{\prime }$ be a random vector with components
\begin{equation*}
C_{i}(n):=\frac{d}{d\lambda }F_{n}^{(i)}(\mathbf{1;}\lambda )\left\vert
_{\lambda =0}\right. =\mathbf{E}_{\mathbf{F}}\varphi _{i}(n),\,i=1,\ldots,m
\end{equation*}%
and let $\mathbf{c}_{n}:=\left( c_{1}(n),\ldots,c_{m}(n)\right) ^{\prime }$
be a random vector with components%
\begin{equation*}
c_{i}(n):=\frac{d}{d\lambda }\phi _{n}^{(i)}(\mathbf{1;}\lambda )\left\vert
_{\lambda =0}\right. =\mathbf{E}_{\mathbf{\phi }}\phi _{i}(n),\,i=1,\ldots,m.
\end{equation*}

Then, by (\ref{GG3}) $C_{m}(n)=c_{m}(n)$ and, for $i<m$%
\begin{equation*}
C_{i}(n)=c_{i}(n)+\sum_{k=i+1}^{m}h_{ik}(n)C_{k}(n).
\end{equation*}%
Hence we get $\ $\ $\mathbf{C}_{n}=\mathbf{c}_{n}+H_{n}^{\Delta }\mathbf{C}%
_{n}$ where
\begin{equation*}
H_{n}^{\Delta }:=\left( h_{ij}(n)I\left( i<j\right) \right) _{i,j=1}^{m}
\end{equation*}%
is the upper triangular matrix generated by $H_{n}$. Thus, $\mathbf{C}%
_{n}=\left( E-H_{n}^{\Delta }\right) ^{-1}\mathbf{c}_{n}$.

The next two statements are easy consequences of Theorem \ref{TfinalZero}
and \ref{ResRep}.

The first theorem gives conditions under which a busy period of a BTPSFPRE
is infinite with positive probability.

\begin{theorem}
\label{TpolSuper}Assume that the MBPRCRE associated with a BTPSFPRE is such
that its underlying MBPRE satisfies conditions of Theorem \ref{Ttan} with $%
\alpha >0$ and, in addition, condition (\ref{Nonzer}) is valid. If $\Phi $
is the total size of the final product accumulated in the BTPSFPRE during a
busy period then $\mathbf{P}\left( \Phi =\infty \right) >0$. In particular,
if the service time of any customer at any station is positive with
probability 1 then the busy period of the BTPSFPRE is infinite with positive
probability.
\end{theorem}

The statement of the results for a BTPSFPRE whose associated MBPRCRE is
subcritical requires more efforts.

\begin{theorem}
\label{Tbusy}Assume that the MBPRCRE associated with a BTPSFPRE is
subcritical, satisfies conditions of Theorem \ref{TailMain} and
\begin{equation*}
\min_{1\leq J\leq m}F_{0}^{(J)}(\mathbf{0;}0)>0.
\end{equation*}%
If the parameter $\kappa $ specified by (\ref%
{DefKappa}) is such that
\begin{equation*}
\max_{1\leq J\leq m}\mathbf{E}\left[ \left( \xi _{J1}+\cdots +\xi
_{Jm}\right) ^{t\vee 1}+\varphi _{J}^{t}(n)\right] <\infty
\end{equation*}%
for any $t\in (0,\kappa ),$ then  there exists a constant $C\in
(0,\infty )$ such that
\begin{equation}
\mathbf{P}\left( \Phi >y\right) \sim Cy^{-\kappa},\,y\rightarrow \infty .
\label{Busy2}
\end{equation}%
In particular, if the final product of any customer is its service time then
the tail distribution of the length $\Phi $ of a busy period of the system
satisfies (\ref{Busy2}).

\begin{corollary}
\label{TpolSub}Under the conditions of Theorem \ref{Tbusy} $\mathbf{E}\Phi
^{x}<\infty $ if and only if $x\in (0,\kappa )$.
\end{corollary}
\end{theorem}

Let us come back to Examples \ref{Examp1} and \ref{Examp2} considered at the
beginning of the paper.

Differentiating (\ref{GateSimple}) at point
$\mathbf{s}=\mathbf{1}$ we see that the matrix $H_{n}$ and the
vector $\mathbf{C}_{n}$ in Example \ref{Examp1} have
elements%
\begin{equation*}
h_{ij}(n)=\gamma _{ij}(n)+\varepsilon _{ij}(n)\mathbf{E}\left[ \tau
_{i}(n)|T_{in}\right] ,\,i,j=1,\ldots ,m
\end{equation*}%
and%
\begin{equation*}
C_{i}(n)=\mathbf{E}\left[ \tau _{i}(n)|T_{in}\right] ,\,i=1,\ldots ,m,
\end{equation*}%
while by differentiating (\ref{Exh1}) at point $\mathbf{s}=\mathbf{1}$ and
taking into account (\ref{ForW}) we conclude after evident transformations
that $h_{ii}(n)=0,i=1,\ldots ,m$ and, for $i\neq j$
\begin{equation}
h_{ij}(n)=\frac{\gamma _{ij}(n)\left( 1-\gamma _{ii}(n)\right) +\varepsilon
_{ij}(n)\mathbf{E}\left[ \tau _{i}(n)|T_{in}\right] }{1-\gamma
_{ii}(n)-\varepsilon _{ii}(n)\mathbf{E}\left[ \tau _{i}(n)|T_{in}\right] },
\label{LL1}
\end{equation}%
if
\begin{equation}
\frac{1-\gamma _{ii}(n)}{\mathbf{E}\left[ \tau _{i}(n)|T_{in}\right] }%
>\varepsilon _{ii}(n),  \label{LL2}
\end{equation}%
and $h_{ij}(n)=\infty ,$ otherwise.

Note that if the quantities $\varepsilon _{ij}(n),\,i,j=1,\ldots ,m$ are
nonrandom, $\gamma _{ii}(n)=0$ with probability 1 for all $i=1,\ldots ,m$,
and $T_{in}$ is an exponential distribution with random parameter $\mu
_{in}, $ relations (\ref{LL1}) and (\ref{LL2}) look as follows:%
\begin{equation*}
h_{ij}(n)=\frac{\mu _{in}\gamma _{ij}(n)+\varepsilon _{ij}(n)}{\mu
_{in}-\varepsilon _{ii}(n)}
\end{equation*}%
and%
\begin{equation*}
\mu _{in}-\varepsilon _{ii}(n)>0.
\end{equation*}%
They are in complete agreement with the respective formulas and restrictions
of Sections 1.1 and 1.2 in \cite{PMPP2008}.

\textbf{Concluding remarks.} Our results give a criterion allowing to answer
the question: when $\Phi $, the total amount of the final product
accumulated in the polling system with zero switchover times during a busy
period has finite or infinite moment of order $x$? In this respect Theorem %
\ref{Tbusy} refines and extends in several directions Theorem 1.1 in \cite%
{PMPP2008}. For instance, we do not require exponentiality of the
service time distributions of customers and prove the mentioned
criterion for a wide
class of polling systems which are not covered by the results of \cite%
{PMPP2008}. Moreover, we even describe the behavior of the tail
distribution of $\Phi $. Unfortunately, such a refinement is
achieved for the expense of transparency of the conditions
involved. The most essential of our hypotheses is Condition $T$,
whose validity is established up to now only for a restricted
class of nonnegative random matrices. The extension of the class
of nonnegative random matrices for which Condition $T$ is valid is
an interesting and challenging problem. Some statements related
with such circle of problems have been obtained quite recently in
a number of papers (see, for instance, \cite{BDGHU08} and
\cite{Gui2006}). Unfortunately, they do not fit the case of
measures concentrated on nonnegative matrices only.

\end{document}